\documentclass[oneside, 12pt]{amsart}

\usepackage{amsfonts}
\usepackage{amssymb}
\usepackage{latexsym}
\usepackage{graphics}
\usepackage[2emode]{psfrag}
\usepackage{mathrsfs}
\usepackage{amsthm}
\usepackage{amsmath}
\usepackage[arrow,matrix,tips,curve,frame]{xy}
\usepackage{enumerate}
\usepackage{hyperref}

\begin{document}

\newcommand{\thlabel}[1]{\label{th:#1}}
\newcommand{\thref}[1]{Theorem~\ref{th:#1}}
\newcommand{\selabel}[1]{\label{se:#1}}
\newcommand{\seref}[1]{Section~\ref{se:#1}}
\newcommand{\lelabel}[1]{\label{le:#1}}
\newcommand{\leref}[1]{Lemma~\ref{le:#1}}
\newcommand{\prlabel}[1]{\label{pr:#1}}
\newcommand{\prref}[1]{Proposition~\ref{pr:#1}}
\newcommand{\colabel}[1]{\label{co:#1}}
\newcommand{\coref}[1]{Corollary~\ref{co:#1}}
\newcommand{\relabel}[1]{\label{re:#1}}
\newcommand{\reref}[1]{Remark~\ref{re:#1}}
\newcommand{\exlabel}[1]{\label{ex:#1}}
\newcommand{\exref}[1]{Example~\ref{ex:#1}}
\newcommand{\delabel}[1]{\label{de:#1}}
\newcommand{\deref}[1]{Definition~\ref{de:#1}}
\newcommand{\eqlabel}[1]{\label{eq:#1}}
\newcommand{\equref}[1]{(\ref{eq:#1})}

\newcommand{\Hom}{{\rm Hom}}
\newcommand{\End}{{\rm End}}
\newcommand{\Ext}{{\rm Ext}}
\newcommand{\Fun}{{\rm Fun}}
\newcommand{\Lax}{{\rm Lax}}
\newcommand{\Cat}{{\bf Cat}}
\newcommand{\Mor}{{\rm Mor}\,}
\newcommand{\Aut}{{\rm Aut}\,}
\newcommand{\Hopf}{{\rm Hopf}\,}
\newcommand{\Ann}{{\rm Ann}\,}
\newcommand{\Ker}{{\rm Ker}\,}
\newcommand{\Coker}{{\rm Coker}\,}
\newcommand{\im}{{\rm Im}\,}
\newcommand{\coim}{{\rm Coim}\,}
\newcommand{\Trace}{{\rm Trace}\,}
\newcommand{\Char}{{\rm Char}\,}
\newcommand{\Mod}{{\bf mod}}
\newcommand{\Spec}{{\rm Spec}\,}
\newcommand{\Span}{{\rm Span}\,}
\newcommand{\sgn}{{\rm sgn}\,}
\newcommand{\Id}{{\rm Id}\,}
\newcommand{\Com}{{\rm Com}\,}
\newcommand{\codim}{{\rm codim}}
\newcommand{\Mat}{{\rm Mat}}
\newcommand{\Coint}{{\rm Coint}}
\newcommand{\Incoint}{{\rm Incoint}}
\newcommand{\can}{{\rm can}}
\newcommand{\sign}{{\rm sign}}
\newcommand{\kar}{{\rm kar}}
\newcommand{\rad}{{\rm rad}}
\newcommand{\Rep}{{\rm Rep}}
\newcommand{\Rmod}{{\sf Rmod}}
\newcommand{\Lmod}{{\sf Lmod}}
\newcommand{\Rcom}{{\sf Rcom}}
\newcommand{\Lcom}{{\sf Lcom}}
\newcommand{\Bicom}{{\sf Bicom}}
\newcommand{\RCOM}{{\sf RCOM}}
\newcommand{\LCOM}{{\sf LCOM}}
\newcommand{\Bim}{{\sf Bim}}
\newcommand{\Bic}{{\sf Bic}}
\newcommand{\Frm}{{\sf Frm}}
\newcommand{\EM}{{\sf EM}}
\newcommand{\REM}{{\sf REM}}
\newcommand{\LEM}{{\sf LEM}}
\newcommand{\fREM}{{\sf fREM}}
\newcommand{\fLEM}{{\sf fLEM}}
\newcommand{\CAT}{{\sf CAT}}
\newcommand{\GAL}{{\sf GAL}}
\newcommand{\Gal}{{\sf Gal}}
\newcommand{\gAL}{{\sf gAL}}
\newcommand{\Tur}{{\sf Tur}}

\def\Ab{\underline{\underline{\rm Ab}}}
\def\lan{\langle}
\def\ran{\rangle}
\def\ot{\otimes}
\def\uot{\ul{\otimes}}
\def\bul{\bullet}
\def\ubul{\ul{\bullet}}

\def\id{\textrm{{\small 1}\normalsize\!\!1}}
\def\To{{\multimap\!\to}}
\def\bigperp{{\LARGE\textrm{$\perp$}}}
\newcommand{\QED}{\hspace{\stretch{1}}
\makebox[0mm][r]{$\Box$}\\}

\def\AA{{\mathbb A}}
\def\BB{{\mathbb B}}
\def\CC{{\mathbb C}}
\def\DD{{\mathbb D}}
\def\EE{{\mathbb E}}
\def\FF{{\mathbb F}}
\def\GG{{\mathbb G}}
\def\HH{{\mathbb H}}
\def\II{{\mathbb I}}
\def\JJ{{\mathbb J}}
\def\KK{{\mathbb K}}
\def\LL{{\mathbb L}}
\def\MM{{\mathbb M}}
\def\NN{{\mathbb N}}
\def\OO{{\mathbb O}}
\def\PP{{\mathbb P}}
\def\QQ{{\mathbb Q}}
\def\RR{{\mathbb R}}
\def\SS{{\mathbb S}}
\def\TT{{\mathbb T}}
\def\UU{{\mathbb U}}
\def\VV{{\mathbb V}}
\def\WW{{\mathbb W}}
\def\XX{{\mathbb X}}
\def\YY{{\mathbb Y}}
\def\ZZ{{\mathbb Z}}

\def\aa{{\mathfrak A}}
\def\bb{{\mathfrak B}}
\def\cc{{\mathfrak C}}
\def\dd{{\mathfrak D}}
\def\ee{{\mathfrak E}}
\def\ff{{\mathfrak F}}
\def\gg{{\mathfrak G}}
\def\hh{{\mathfrak H}}
\def\ii{{\mathfrak I}}
\def\jj{{\mathfrak J}}
\def\kk{{\mathfrak K}}
\def\ll{{\mathfrak L}}
\def\mm{{\mathfrak M}}
\def\nn{{\mathfrak N}}
\def\oo{{\mathfrak O}}
\def\pp{{\mathfrak P}}
\def\qq{{\mathfrak Q}}
\def\rr{{\mathfrak R}}
\def\ss{{\mathfrak S}}
\def\tt{{\mathfrak T}}
\def\uu{{\mathfrak U}}
\def\vv{{\mathfrak V}}
\def\ww{{\mathfrak W}}
\def\xx{{\mathfrak X}}
\def\yy{{\mathfrak Y}}
\def\zz{{\mathfrak Z}}

\def\aaa{{\mathfrak a}}
\def\bbb{{\mathfrak b}}
\def\ccc{{\mathfrak c}}
\def\ddd{{\mathfrak d}}
\def\eee{{\mathfrak e}}
\def\fff{{\mathfrak f}}
\def\ggg{{\mathfrak g}}
\def\hhh{{\mathfrak h}}
\def\iii{{\mathfrak i}}
\def\jjj{{\mathfrak j}}
\def\kkk{{\mathfrak k}}
\def\lll{{\mathfrak l}}
\def\mmm{{\mathfrak m}}
\def\nnn{{\mathfrak n}}
\def\ooo{{\mathfrak o}}
\def\ppp{{\mathfrak p}}
\def\qqq{{\mathfrak q}}
\def\rrr{{\mathfrak r}}
\def\sss{{\mathfrak s}}
\def\ttt{{\mathfrak t}}
\def\uuu{{\mathfrak u}}
\def\vvv{{\mathfrak v}}
\def\www{{\mathfrak w}}
\def\xxx{{\mathfrak x}}
\def\yyy{{\mathfrak y}}
\def\zzz{{\mathfrak z}}

\newcommand{\aA}{\mathscr{A}}
\newcommand{\bB}{\mathscr{B}}
\newcommand{\cC}{\mathscr{C}}
\newcommand{\dD}{\mathscr{D}}
\newcommand{\eE}{\mathscr{E}}
\newcommand{\fF}{\mathscr{F}}
\newcommand{\gG}{\mathscr{G}}
\newcommand{\hH}{\mathscr{H}}
\newcommand{\iI}{\mathscr{I}}
\newcommand{\jJ}{\mathscr{J}}
\newcommand{\kK}{\mathscr{K}}
\newcommand{\lL}{\mathscr{L}}
\newcommand{\mM}{\mathscr{M}}
\newcommand{\nN}{\mathscr{N}}
\newcommand{\oO}{\mathscr{O}}
\newcommand{\pP}{\mathscr{P}}
\newcommand{\qQ}{\mathscr{Q}}
\newcommand{\rR}{\mathscr{R}}
\newcommand{\sS}{\mathscr{S}}
\newcommand{\tT}{\mathscr{T}}
\newcommand{\uU}{\mathscr{U}}
\newcommand{\vV}{\mathscr{V}}
\newcommand{\wW}{\mathscr{W}}
\newcommand{\xX}{\mathscr{X}}
\newcommand{\yY}{\mathscr{Y}}
\newcommand{\zZ}{\mathscr{Z}}

\newcommand{\Aa}{\mathcal{A}}
\newcommand{\Bb}{\mathcal{B}}
\newcommand{\Cc}{\mathcal{C}}
\newcommand{\Dd}{\mathcal{D}}
\newcommand{\Ee}{\mathcal{E}}
\newcommand{\Ff}{\mathcal{F}}
\newcommand{\Gg}{\mathcal{G}}
\newcommand{\Hh}{\mathcal{H}}
\newcommand{\Ii}{\mathcal{I}}
\newcommand{\Jj}{\mathcal{J}}
\newcommand{\Kk}{\mathcal{K}}
\newcommand{\Ll}{\mathcal{L}}
\newcommand{\Mm}{\mathcal{M}}
\newcommand{\Nn}{\mathcal{N}}
\newcommand{\Oo}{\mathcal{O}}
\newcommand{\Pp}{\mathcal{P}}
\newcommand{\Qq}{\mathcal{Q}}
\newcommand{\Rr}{\mathcal{R}}
\newcommand{\Ss}{\mathcal{S}}
\newcommand{\Tt}{\mathcal{T}}
\newcommand{\Uu}{\mathcal{U}}
\newcommand{\Vv}{\mathcal{V}}
\newcommand{\Ww}{\mathcal{W}}
\newcommand{\Xx}{\mathcal{X}}
\newcommand{\Yy}{\mathcal{Y}}
\newcommand{\Zz}{\mathcal{Z}}

\def\units{{\mathbb G}_m}
\def\rightact{\hbox{$\leftharpoonup$}}
\def\leftact{\hbox{$\rightharpoonup$}}

\def\*C{{}^*\hspace*{-1pt}{\Cc}}

\def\text#1{{\rm {\rm #1}}}

\def\smashco{\mathrel>\joinrel\mathrel\triangleleft}
\def\cosmash{\mathrel\triangleright\joinrel\mathrel<}

\def\ol{\overline}
\def\ul{\underline}
\def\dul#1{\underline{\underline{#1}}}
\def\Nat{\dul{\rm Nat}}
\def\Set{\dul{\rm Set}}

\renewcommand{\subjclassname}{\textup{2000} Mathematics Subject
     Classification}

\newtheorem{proposition}{Proposition}[section]
\newtheorem{lemma}[proposition]{Lemma}
\newtheorem{corollary}[proposition]{Corollary}
\newtheorem{theorem}[proposition]{Theorem}

\theoremstyle{definition}
\newtheorem{Definition}[proposition]{Definition}
\newtheorem{example}[proposition]{Example}
\newtheorem{examples}[proposition]{Examples}

\theoremstyle{remark}
\newtheorem{remarks}[proposition]{Remarks}
\newtheorem{remark}[proposition]{Remark}

\title{Galois Theory in Bicategories}
\date{\today}
\author{J. G\'omez-Torrecillas}
\address{Departamento de \'{A}lgebra, Universidad
de Granada, E-18071 Granada, Spain} \email{gomezj@ugr.es}
\urladdr{www.ugr.es/\~{}gomezj}

\author{J. Vercruysse}
\address{Faculty of Engineering, Vrije Universiteit Brussel (VUB), B-1050 Brussels, Belgium}
\email{jvercruy@vub.ac.be}
\urladdr{homepages.vub.ac.be/\~{}jvercruy}
\thanks{The research of the first author is supported by the grants MTM2007-61673 from the Ministerio
de Edu\-ca\-ci{\'o}n y Ciencia of Spain, and P06-FQM-01889 from the
Consejer{\'\i}a de Innovaci{\'o}n, Ciencia y Empresa of Andaluc{\'\i}a, with
funds from FEDER (Uni{\'o}n Europea)}
\thanks{The second author is Postdoctoral Fellow of the Fund for Scientific Research--Flanders
(Belgium) (F.W.O.--Vlaanderen).}


\begin{abstract}
We develop a Galois (descent) theory for comonads within the
framework of bicategories. We give generalizations of Beck's theorem
and the Joyal-Tierney theorem. Many examples are provided, including
classical descent theory, Hopf-Galois theory over Hopf algebras and
Hopf algebroids, Galois theory for corings and group-corings, and
Morita-Takeuchi theory for corings. As an application we construct a
new type of comatrix corings based on (dual) quasi bialgebras.
\end{abstract}

\maketitle

\section*{Introduction}

The classical Galois Theory on field extensions has been
generalized in many directions. For instance, it has been extended
to a Galois theory for commutative rings by Auslander and Goldman
\cite{AusGold} and by Chase, Harrison and Rosenberg \cite{CHR}. A
group action can be generalized to a Hopf algebra (co)action. This
leads to the Hopf-Galois theory, developed first for finitely
generated and projective Hopf algebras (see \cite{CS} and
\cite{KT}) and later for arbitrary Hopf algebras (see \cite{DT}
and \cite{Schneider90}). During the nineties, the theory of Hopf
algebras went through a range of generalizations, such as
Doi-Koppinen structures \cite{Doi:DK}, \cite{Kop:DK} and entwining
structures \cite{BrzezinskiH} to arrive at the theory of corings
and comodules \cite{Sw}, which provides a general framework to
explain many results of Hopf algebra theory in a simple and
clarifying way. In this respect, it is no surprise that
Hopf-Galois theory has a formulation in terms of corings. This was
shown in \cite{Brz:structure}, where a Galois theory is developed
for corings with a grouplike element. To a ring morphism $i:B\to
A$, we can associate an $A$--coring, the so called canonical
Sweedler coring. A morphism from this coring to another $A$--coring
$\cc$ is completely determined by a grouplike element $g\in\cc$.
When this morphism is an isomorphism we say that $(\cc,g)$ is a
Galois coring. We can construct a pair of adjoint functors between
the categories $\Mm_B$ and $\Mm^\cc$ and formulate sufficient and
necessary conditions for this pair to be an equivalence of
categories. El Kaoutit and the first author \cite{EGT:comatrix}
introduced a yet more general version of Galois theory, replacing
the grouplike element by a right $\cc$--comodule $\Sigma$ which is
finitely generated and projective as a right $A$--module. This
generalizes the Galois theory of \cite{Brz:structure}, and, in
particular, allows to characterize those corings whose category of
comodules has a finitely generated projective generator. To this
end, an $A$--coring, called the comatrix coring, is constructed out
of $\Sigma$. We call $\Sigma$ a Galois comodule if the canonical
homomorphism from the comatrix coring to $\cc$ is bijective (see
also \cite{Brz:galcom}, \cite{CDV}). In general, there exists an
adjunction between the categories $\Mm_T$ and $\Mm^\cc$, where
$T=\End^\cc(\Sigma)$, and this adjunction induces an equivalence of
categories if and only if $\cc$ is Galois and ${}_T\Sigma$ is
faithfully flat \cite[Theorem 3.2]{EGT:comatrix}. Several attempts
were made to drop or weaken the finiteness condition on $\Sigma$
and construct infinite versions of Galois theory for comodules. El
Kaoutit and the first author introduced Galois comodules that are
(infinite) direct sums of finitely generated and projective right
$A$--modules \cite{EGT:infinite}. Caenepeel, De Groot and the
second author generalized this to a method to construct corings
out of colimits. Both authors developed a theory of Galois
comodules over firm rings \cite{GTV}. Wisbauer introduced a
functorial definition for a Galois comodule \cite{Wis:galcom}.

In a slightly different direction, there has been a similar evolution in (categorical) descent
theory. Both theories are closely related (this connection is discussed in \cite{GT:comonad}),
and some authors would prefer to term the theory of this paper `descent theory'. However, we
reserve the name `descent theory' for the special situation that is treated in
\seref{descent}. Our motivation to keep the term Galois theory is founded by the evolution in
Hopf-Galois theory described above. Descent theory investigates the extension of scalars
functor $-\ot_BA:\Mm_B\to\Mm_A$ associated to a homomorphism of rings $B\to A$, in particular
it looks for sufficient and necessary conditions for this functor to be comonadic
(cotripleable). A first important theorem in this respect, is Beck's theorem that gives
sufficient and necessary conditions for a functor with a right adjoint to be comonadic (see
e.g. \cite{BarrWells:ttt}, \cite{McLane}). Another interesting and more general result states
that the extension of scalars functor for commutative rings $A$ and $B$ is comonadic if and
only if $B\to A$ is a pure monomorphism of $B$--modules. Although this theorem is presently
known as the theorem of Joyal-Tierney, it was never published as such; a proof of this theorem
can be found in \cite{Mes:pure}. During the last few years, the Joyal-Tierney theorem is
generalized in several ways. In \cite{Cae:descent}, a non-commutative version is presented, in
\cite{CDV} a generalization is formulated where the ring extension is replaced by a $B$-$A$--bimodule, where $A$ is non-commutative and $B$ is commutative, in \cite{Mes:extension} $B$ is
allowed to be a (non-commutative) separable algebra. In \cite{JanTho:descIII} a categorical
version of the theorem is presented.

In this paper we propose a Galois theory for comonads in the
general setting of bicategories. Our approach rests upon the set
up proposed in \cite{LacStr} in the light of the theory developed
in \cite{BrzElkGT} for bimodules and corings. Our work intends to
provide not only a transparent view on the interactions between
the different approaches in the recent development on Galois
comodules, but we also hope to shed new light on the relationship
between the coring theory and the theory of comonads
\cite{GT:comonad}. Moreover several other versions and
generalizations of Hopf-Galois theory, which have been formulated
during the last years (such as equivalences between categories of
comodules \cite{BrzElkGT}, \cite{Zarouali:phd}, Galois theory for
$C$--rings \cite{BrzTur:Crings} and Galois theory for group-corings
\cite{CJW}) fit perfectly within our general framework.

One final remark on notation: in any category $\Cc$, we will denote the identity morphism on an object $X\in\Cc$ again by $X$.

\section{Elementary definitions and notation}\selabel{definitions}

Recall from \cite{Ben} that a \emph{bicategory} $\Bb$ consists of
the following data.
\begin{enumerate}[(i)]
\item A class of objects $A,B,\ldots$ which are called $0$--cells
(or objects). \item For every two objects $A$ and $B$, there
exists a category $\Hom_\Bb(A,B)=\Hom(A,B)$, whose class of
objects we denote by $\Hom_1(A,B)$ and which are called $1$--cells.
We denote $f:A\to B$ for an $1$--cell $f\in \Hom_1(A,B)$. Take two
$1$--cells $f, g\in\Hom_1(A,B)$. The set of morphisms from $f$ to
$g$ in the category $\Hom(A,B)$ is denoted by ${^A\Hom_2^B}(f,g)$.
We call these morphisms $2$--cells and denote them as $\alpha:f\to
g$. We will denote the composition of morphisms in the category
$\Hom(A,B)$ by $\circ$, i.e. for all $f,g,h\in\Hom_1(A,B)$ such
that $\alpha:f\to g$ and $\beta:g\to h$, we have $\beta\circ
\alpha:f\to h$. This composition will now be called the
\emph{vertical} composition of $2$--cells. \item For any three
objects $A,B, C\in \Bb$, there exists a functor
$$c_{ABC}:\Hom(A,B) \times \Hom(B,C) \to \Hom(A,C).$$
For all $f\in\Hom_1(A,B)$ and $g\in\Hom_1(B,C)$, we denote
$c_{ABC}(f,g)=f\bul_B g\in\Hom_1(A,C)$. For all
$\alpha\in{^A\Hom_2^B}(f,g)$ and $\beta\in{^B\Hom_2^C}(h,k)$, we
denote $c_{ABC}(\alpha,\beta)=\alpha \ubul_B\beta:f\bul_B h\to
g\bul_B k$. This composition will be called the \emph{horizontal}
composition of $2$--cells. \item For any object $A\in\Bb$, there
exists a functor
$$\id_A : {\bf 1} \to \Hom(A,A),$$
where ${\bf 1}$ denotes the discrete category with one object $*$.
We will denote $\id_A(*)$ just by $\id_A$. \item For any four
objects $A,B,C,D\in\Bb$, there exists a natural isomorphism
\begin{equation}\eqlabel{associativity}
\alpha_{ABCD}:c_{ACD}\circ(c_{ABC}\times \Hom(C,D))\to
c_{ABD}\circ(\Hom(A,B)\times c_{BCD}).
\end{equation}
For any two objects $A,B\in \Bb$, there exist two natural
isomorphisms
\begin{eqnarray}\eqlabel{rho}
\lambda_{AB}:\Hom(A,B)&\to& c_{AAB}\circ (\id_A\times
\Hom(A,B)),\\ \rho_{AB}:\Hom(A,B)&\to&
c_{ABB}\circ(\Hom(A,B)\times \id_B).\nonumber
\end{eqnarray}
\end{enumerate}
All these data are required to satisfy some compatibility
(associativity and coherence) conditions. We refer to e.g.\
\cite{Ben},
where the notion of a bicategory was introduced, or \cite[section 7.7]{Bo}.\\

For all objects $A,B,C\in\Bb$, we obtain from the functorality of
$c_{ABC}$ the \emph{interchange law}, i.e.
\begin{equation}\eqlabel{interchange}
(\alpha\ubul_B\beta)\circ(\gamma\ubul_B\delta)=(\alpha\circ\gamma)\ubul_B(\beta\circ\delta),
\end{equation}
for $\alpha\in{^A\Hom_2^B}(a,c)$, $\beta\in{^B\Hom_2^C}(b,d)$,
$\gamma\in{^A\Hom_2^B}(c,e)$ and $\delta\in{^B\Hom_2^C}(d,f)$.
From \equref{interchange} one immediately deduces that for all
$\alpha\in{^A\Hom_2^B}(a,c)$ and $\beta\in{^B\Hom_2^C}(b,d)$,
\begin{equation}\eqlabel{interchange2}
(\alpha\ubul_B b)\circ(c\ubul_B
\beta)=\alpha\ubul\beta=(a\ubul_B\beta)\circ(\alpha\ubul_B d).
\end{equation}

A \emph{$2$--category} is a bicategory such that the isomorphisms
$\alpha_{ABCD}$, $\lambda_{AB}$ and $\rho_{AB}$ are identities for
all choices of $A,B,C,D$. In particular, $\id_A=A$ for all objects
$A$ of a $2$--category.

To any bicategory $\Bb$ one can associate new bicategories denoted
by $\Bb^{op}$, $\Bb^{co}$ and $\Bb^{coop}$. These are constructed
by taking respectively opposite composition for the $1$--cells, for
vertical composition of $2$--cells and for both.

Recall that
a \emph{comonad} $\cc=(A,c,\Delta_c,\varepsilon_c)$ in $\Bb$
consists of a $0$--cell $A$, an $1$--cell $c\in\Hom_1(A,A)$ and two
$2$--cells $\Delta_c\in{^A\Hom_2^A}(c,c\bullet_Ac)$ and
$\varepsilon_c\in{^A\Hom_2^A}(c,\id_A)$ such that the diagrams
\[
\xymatrix{
& c\bullet_Ac \ar[r]^-{\Delta_c\ubul_Ac}& (c\bullet_Ac)\bullet_Ac \ar[dd]^{\alpha_{AAAA}}\\
c \ar[ru]^-{\Delta_c}
\ar[rd]_{\Delta_c} \\
& c\bullet_Ac \ar[r]_-{c\ubul_A\Delta_c} & c\bullet_A(c\bullet_Ac)
} \qquad \xymatrix{
c \ar[rrd]^-{\Delta_c} \ar[rr]^{\rho_{AA}} \ar[d]_{\lambda_{AA}} && c\bullet_A\id_A \\
\id_A\bullet_Ac && c\bullet_Ac \ar[ll]^{\varepsilon_c\ubul_Ac}
\ar[u]_{c\ubul_A\varepsilon_c} }
\]
are commutative in $\Hom(A,A)$. For each object $\Omega$ of $\Bb$,
the comonad $\cc$ in $\Bb$ induces a comonad on the category
$\Hom(\Omega,A)$ (in the sense of \cite[Definition 3.1]{BarrWells:ttt})
\[
\delta : G \rightarrow G^2, \qquad \epsilon : G \rightarrow 1,
\]
where
\[
G = - \bullet_A c, \quad \delta = \alpha^{-1}_{\Omega AAA} \circ
(- \ubul_A \Delta_c), \quad \text{and } \quad \epsilon =
\rho_{\Omega A}^{-1} \circ (- \ubul_A \varepsilon_c).
\]
Let $\Rcom(\Omega,\cc)$ denote the Eilenberg-Moore category for this
comonad. Its objects will be called \emph{right $\cc$--comodules of
$\Omega$--type} or \emph{$\Omega$-$\cc$--comodules}, and consist of
couples $(m,\rho^m)$ where $m\in\Hom_1(\Omega,A)$ for some $0$--cell
$\Omega$ in $\Bb$ and $\rho^m\in{^\Omega\Hom_2^A}(m,m\bullet_Ac)$
such that $\alpha_{\Omega
AAA}\circ(\rho^m\ubul_Ac)\circ\rho^m=(m\ubul_A\Delta_c)\circ\rho^m$
and $(m\ubul_A\varepsilon_c)\circ\rho^m=\rho_{\Omega A}$. A morphism
of $\Omega$-$\cc$--comodules, $\psi:(m,\rho^m)\to(n,\rho^n)$
consists of a $2$--cell $\psi\in{^\Omega\Hom_2^A}(m,n)$ such that
$\rho^n\circ\psi=(\psi\ubul_Ac)\circ\rho^m$. We will denote the set
of all morphisms of $\Omega$-$\cc$--comodules between $m$ and $n$ by
$^\Omega\Hom^\cc(m,n)$. Analogously, our comonad $\cc$ induces, for
each $0$--cell $\Omega$, a comonad on $\Hom(A,\Omega)$, whose
corresponding Eilenberg-Moore category will be denoted by
$\Lcom(\cc,\Omega)$, and their objects will be referred to as
\emph{left $\cc$--comodules of $\Omega$--type} or
\emph{$\cc$-$\Omega$--comodules}.

Given comonads $\cc = (A,c,\Delta_c,\varepsilon_c)$ and $\dd =
(B,d,\Delta_d,\varepsilon_d)$, we define a
$\cc$-$\dd$--\emph{bicomodule} as a three-tuple
$(m,\lambda^m,\rho^m)$ where $(m,\lambda^m)$ is a
$\cc$-$B$--comodule and $(m,\rho^m)$ is an $A$-$\dd$--comodule, such
that
$\alpha_{AABB}\circ(\lambda^m\ubul_Bd)\circ\rho^m=(c\ubul_A\rho^m)\circ\lambda^m$.
These bicomodules are the objects of a category $\Bicom(\cc,\dd)$,
whose morphisms are the $2$--cells that are both morphisms of left
$\cc$--comodules and of right $\dd$--comodules.

A \emph{monad} and \emph{modules} over a monad in a bicategory $\Bb$ can be defined as a comonad and comodules over a comonad in the co-opposite bicategory $\Bb^{co}$, which consists of the same data as $\Bb$, except that one considers the opposite composition of the vertical composition of $2$--cells in $\Bb$.

A pseudo-functor $F : \Bb \rightarrow \Cc$ between bicategories
$\Bb$ and $\Cc$ assigns a $0$--cell $FA$ of $\Cc$ to each
$0$--cell $A$ of $\Bb$, and a functor
$$F_{AB}:\Hom_\Bb(A,B) \to \Hom_\Cc(FA,FB);$$
for every pair of $0$--cells $A,B$ of $\Bb$. The pseudo-functor
preserves the horizontal composition only up to isomorphisms, in
the sense that for every three-tuple of $0$--cells $A,B,C$ in
$\Bb$, there exist natural isomorphisms
\begin{equation}\eqlabel{laxfunctor}
\gamma_{ABC}:c_{FA,FB,FC}\circ(F_{AB}\times F_{BC})\to F_{AC}\circ
c_{ABC};
\end{equation}
\begin{equation}\eqlabel{laxfunctor2}
\xi_A:\id_{FA}\to F_{AA}\circ \id_A;
\end{equation}
subject to suitable associativity and coherence axioms (see e.g.
\cite[Section 7.5]{Bo}). If $P$ is a property of functors, then we say that a pseudo functor $F$ satisfies the \emph{local} $P$--property if and only if $F_{AB}$ satisfies the propery $P$ for all choices of $A$ and $B$ (e.g. $F$ is locally faithful, locally exact).

Given a comonad $\cc =
(A,c,\Delta,\varepsilon)$ in $\Bb$, the pseudo-functor $F:\Bb\to\Cc$ induces a
comonad in $\Cc$ given by $$F\cc = (FA,Fc,\gamma_{AAA}^{-1} \circ
F\Delta, \xi_A^{-1} \circ F\varepsilon),$$ and, for each object
$\Omega$ of $\Bb$, a functor $F_{\cc} : \Rcom(\Omega,\cc)
\rightarrow \Rcom(F\Omega,F\cc)$. If all $F_{AB}$'s are
equivalences (e.g. if $F$ is a biequivalence of bicategories),
then $F_{\cc}$ is an equivalence for every comonad $\cc$. Since
every bicategory is biequivalent to a $2$--category (see
\cite{Str:book}, \cite{McLanePar}), we can, without loss of
generality, restrict ourselves to the case of comonads in
$2$--categories when studying categories of comodules. In fact,
our former argument could have been replaced by the more general
``Coherence Theorem'', which asserts that all diagrams that are
constructed out of the associativity and identity isomorphisms
commute in any bicategory. Those readers familiar enough with the
Coherence Theorem may consider that we will rely on the
$2$--categorical calculus when we are dealing with bicategories.

\section{The 2--categories of Comonads}\selabel{bicatcomonad}

Throughout this section, $\Bb$ will denote a $2$--category.
Following \cite{Str} and \cite{LacStr}, we will define several
$2$--categories whose $0$--cells are comonads in $\Bb$. Their
$1$--cells will be comonad-morphisms, in the sense of Definition
\ref{comonadmorphism}. They will encode certain functors between
categories of comodules (see \seref{Galois}). In this section, we
show that the comonad-morphisms can be understood as bicomodules
(\leref{bicomodule}).

\begin{Definition}\label{comonadmorphism}
Let $\cc=(A,c,\Delta_c,\varepsilon_c)$ and
$\dd=(B,d,\Delta_d,\varepsilon_d)$ be two comonads in $\Bb$. A
\emph{right comonad-morphism} from $\dd$ to $\cc$ is a pair
$(q,\alpha)$, consisting of $q\in\Hom_1(B,A)$ and
$\alpha\in{^B\Hom_2^A}(d\bul_Bq,q\bul_Ac)$ such that the following
diagrams commute
\[
\xymatrix{
d\bul_Bq \ar[rr]^-{\alpha} \ar[d]_{\Delta_d\ubul_Bq} && q\bul_Ac \ar[d]^{q\ubul_A\Delta_c}\\
d\bul_Bd\bul_Bq \ar[r]_{d\ubul_B\alpha} & d\bul_Bq\bul_Ac
\ar[r]_{\alpha\ubul_Ac} & q\bul_Ac\bul_Ac }\quad \xymatrix{
d\bul_Bq \ar[dr]_{\varepsilon_d\ubul_Bq} \ar[rr]^-{\alpha} && q\bul_Ac \ar[dl]^{q\ubul_A\varepsilon_c}\\
& q }
\]
In a symmetric way, a \emph{left comonad-morphism} from $\dd$ to
$\cc$ is a pair $(p,\beta)$, where $p\in\Hom_1(A,B)$ and
$\beta\in{^A\Hom_2^B}(p\bul_Bd,c\bul_Ap)$, such that the following
diagrams commute
\[
\xymatrix{
p\bul_Bd \ar[rr]^-{\beta} \ar[d]_{p\ubul_B\Delta_d} && c\bul_Ap \ar[d]^{\Delta_c\ubul_Ap}\\
p\bul_Bd\bul_Bd \ar[r]_{\beta\ubul_Bp} & c\bul_Ap\bul_Bd
\ar[r]_{c\ubul_A\beta} & c\bul_Ac\bul_Ap }\quad \xymatrix{
p\bul_Bd \ar[dr]_{p\ubul_B\varepsilon_d} \ar[rr]^-{\beta} && c\bul_Ap \ar[dl]^{\varepsilon_c\ubul_Ap}\\
& p }
\]
\end{Definition}

\begin{examples}\exlabel{trivcomonad}
\begin{enumerate}[(i)]
\item Let
$\varphi:\cc=(A,c,\Delta,\varepsilon)\to\cc'=(A,c',\Delta',\varepsilon')$
be a \emph{homomorphism of $A$--comonads}, i. e., a $2$--cell
$\varphi : c \rightarrow c'$ such that
$\Delta_{c'}\circ\varphi=(\varphi\ubul_A\varphi)\circ\Delta_c$ and
$\varepsilon_{c'}\circ\varphi=\varepsilon_c$. Then $(A,\varphi)$
is both a left and a right comonad-morphism from $\cc$ to $\cc'$.
\item From the definition one can immediately see that a right
comonad-morphism $(m,\rho)$ from $\bb$ to $\cc$, where
$\bb=(B,B,B,B)$ is the trivial $B$--comonad, is nothing else than a
right $\cc$--comodule of $B$--type. Similarly, left comonad-morphisms
from $\bb$ to $\cc$ are exactly left $\cc$--comodules.
More generally, we have the following result.
\end{enumerate}
\end{examples}

\begin{lemma}\lelabel{bicomodule}
Let $\cc=(A,c,\Delta_c,\varepsilon_c)$ and $\dd=(B,d,\Delta_d,\varepsilon_d)$ be comonads in $\Bb$. Then the following statements hold
\begin{enumerate}[(i)]
\item $(q,\alpha)$ is a right comonad-morphism from $\dd$ to $\cc$
if and only if $d\bul_Bq$ is a $\dd\hbox{-}\cc$--bicomodule. \item
$(p,\beta)$ is a left comonad-morphism from $\dd$ to $\cc$ if and
only if $p\bul_Bd$ is a $\cc\hbox{-}\dd$--bicomodule.
\end{enumerate}
\end{lemma}

\begin{proof}
Suppose first that $(q,\alpha)$ and $(p,\beta)$ are comonad-morphisms. The coactions on $d\bul_Bq$ and $p\bul_Bd$ are given by the following formulas,
\begin{eqnarray*}
\rho^{d\bul_Bq}:&& \xymatrix{d\bul_Bq \ar[r]^-{\Delta_d\ubul_Bq} & d\bul_Bd\bul_Bq \ar[r]^-{d\ubul_B\alpha} & d\bul_Bq\bul_Ac};\\
\lambda^{d\bul_Bq}:&&\xymatrix{d\bul_Bq \ar[r]^-{\Delta_d\ubul_Bq} & d\bul_Bd\bul_Bq};\\
\rho^{p\bul_Bd}:&&\xymatrix{p\bul_Bd\ar[r]^-{p\ubul_B\Delta_d} & p\bul_Bd\bul_Bd};\\
\lambda^{p\bul_Bd}:&&\xymatrix{p\bul_Bd\ar[r]^-{p\ubul_B\Delta_d} & p\bul_Bd\bul_Bd \ar[r]^-{\beta\ubul_Bd} & c\bul_Ap\bul_Bd}.
\end{eqnarray*}
We will only check the coassociativity of $\rho^{d\bul_Bq}$ and
leave other verifications (the coassociativity and counit
conditions, as well as the compatibility between left and right
coaction) to the reader, since they are all similar computations.
Consider the following diagram, of which the outer quadrangle
expresses the coassociativity of $\rho^{d\bul_Bq}$.
\[
\xymatrix{
d\bul_Bq \ar[rr]^{\Delta_d\ubul_Bq} \ar[d]^{\Delta_d\ubul_Bq} && d\bul_Bd\bul_Bq \ar[rr]^{d\ubul_B\alpha} &&
d\bul_Bq\bul_Ac \ar[d]^{\Delta_d\ubul_Bq\ubul_Ac}\\
d\bul_Bd\bul_Bq \ar[rr]^-{\Delta_d\ubul_Bd\ubul_Bq}_-{d\ubul_B\Delta_d\ubul_Bq} \ar[d]_{d\ubul_B\alpha}
&& d\bul_Bd\bul_Bd\bul_Bq \ar[rr]^{d\ubul_Bd\ubul_B\alpha} && d\bul_Bd\bul_Bq\bul_Ac \ar[d]^{d\ubul_B\alpha\ubul_Ac}\\
d\bul_Bq\bul_Ac \ar[rrrr]_{d\ubul_Bq\ubul_A\Delta_c} &&&& d\bul_Bq\bul_Ac\bul_Ac
}
\]
The lower quadrangle of this diagram commutes by the definition of
a right comonad-morphism, and the upper quadrangle by application
of \equref{interchange2}.

Now suppose that $d\bul_B q$ is a $\dd$-$\cc$--bicomodule. We will
prove that $(q,\alpha)$ is a comonad-morphism, where
\[
\xymatrix{
\alpha:d\bul_B q \ar[rr]^-{d\ubul_B\rho} && d\bul_B q\bul_A c \ar[rr]^-{\varepsilon_d\ubul_Bq\ubul_Ac} && q\bul_A c.
}
\]
We have to check the commutativity of the outer quadrangle of the following diagram.
\[
\xymatrix{
d\bul_B q \ar[dd]_-{d\ubul_B\rho} \ar[r]^-{\Delta_d\ubul_Bq} & d\bul_Bd\bul_Bq \ar[d]^-{d\ubul_B\rho}\\
& d\bul_Bd\bul_Bq\bul_Ac \ar[d]^-{d\ubul_B\varepsilon_d\ubul_Bq\ubul_Ac} \\
d\bul_B q\bul_A c \ar[ru]^-{\Delta_d\ubul_Bq\ubul_Ac} \ar@{=}[r] \ar[rd]_-{d\ubul_Bq\ubul_A\Delta_c} \ar[dd]_-{\varepsilon_d\ubul_Bq\ubul_Ac}
& d\bul_Bq\bul_Ac \ar[d]^-{d\ubul\rho}  \\
&  d\bul_Bq\bul_Ac\bul_Ac \ar[d]^-{\varepsilon_d\ubul_Bq\ubul_Ac\ubul_Ac}\\
q\bul_A c  \ar[r]_-{q\ubul_A\Delta_c}
 & q\bul_Ac\bul_Ac
}
\]
The upper inner quadrangle commutes since $d\bul_Bq$ is a bicomodule. The upper triangle commutes by the counital property of the comonad $\dd$. The lower triangle commutes on the image of the incoming arrow $d\ubul_B\rho$ on the left as an application of the coassociativity condition on $d\bul_Bq$ as a right $\cc$--comodule. Finally, applying \equref{interchange2}, we find that the lower quadrangle commutes.

The proof for $p\bul_B d$ is similar.
\end{proof}

The $2$--category of comonads in $\Bb$ introduced in \cite{Str} is
the following:

\begin{Definition}By $\RCOM(\Bb)$ we denote the \emph{right $2$--category of comonads} in $\Bb$.
This $2$--category is defined as follows.
\begin{enumerate}[(i)]
\item A $0$--cell in $\RCOM(\Bb)$ is a comonad in $\Bb$; \item a
$1$--cell in $\RCOM(\Bb)$ is a right comonad-morphism; \item a
$2$--cell in $\RCOM(\Bb)$ between two right comonad-morphisms
$(q,\alpha)$ and $(q',\alpha')$ from $\dd$ to $\cc$ is a $2$--cell
$\sigma:q\to q'$ in $\Bb$ that makes the following diagram commute
\[
\xymatrix{
d\bul_Bq \ar[r]^{d\ubul_B\sigma} \ar[d]_{\alpha} & d\bul_Bq' \ar[d]^{\alpha'}\\
q\bul_Ac \ar[r]_{\sigma\ubul_Ac} & q'\bul_Ac ;}
\]
\item the composition of one-cells is defined as follows. Consider comonads $\cc=(A,c,\Delta_c,\varepsilon_c)$,
$\dd=(B,d,\Delta_D,\varepsilon_d)$ and
$\ee=(C,e,\Delta_e,\varepsilon_e)$. Let $(q,\alpha): \cc\to\dd$
and $(q', \alpha'):\dd\to \ee$ be two comonad-morphisms. Then we
define a new comonad-morphism $(q,\alpha)\bul_\dd(q',\alpha') =
(q'\bul_Bq,(q'\ubul_A\alpha) \circ(\alpha'\ubul_B q)):\cc\to\ee$.
\end{enumerate}
Considering left comonad-morphisms, one obtains the left
$2$--category of comonads $\LCOM(\Bb)$.
\end{Definition}

In \cite{LacStr}, an alternative $2$--category of comonads is
proposed, which contains the same $0$--cells and $1$--cells, but
different $2$--cells. In \cite[p. 249]{LacStr} there were given two
equivalent descriptions of the $2$--cells, called the reduced and
the unreduced form. In \cite[Proposition 2.2]{BrzElkGT}, where the
dual case for $\Bb=\Bim$, the bicategory of bimodules, was
considered, the unreduced form is being interpreted as a morphism
of bicomodules. We think that the treatment of $2$--cells as
bicomodules correspondes to yet a different description of the
$2$--cells than the ones that can be found in \cite{LacStr}, as the
following lemma explains.

\begin{lemma}\lelabel{sigma}
Let $(q,\alpha)$ and $(q',\alpha')$ be two right comonad-morphisms
from $\dd$ to $\cc$ in $\Bb$. There exists a bijective
correspondence between the following objects:
\begin{enumerate}[(i)]
\item a $2$--cell $\sigma:d\bul_Bq\to q'$ in $\Bb$ making the following diagram commutative
\[
\xymatrix{
d\bul_Bq \ar[r]^-{\Delta_d\ubul_Bq} \ar[d]_{\Delta_d\ubul_Bq} &d\bul_Bd\bul_Bq \ar[r]^-{d\ubul_B\sigma} &d\bul_Bq' \ar[d]^{\alpha'}\\
d\bul_Bd\bul_Bq \ar[r]_-{d\ubul_B\alpha} & d\bul_Bq\bul_Ac \ar[r]_-{\sigma\ubul_Ac} & q'\bul_Ac
}
\]
\item a $2$--cell $\widetilde{\sigma}:d\bul_Bq\to q'\bul_Ac$ in $\Bb$ making the following two diagrams commutative
\[
\xymatrix{
d\bul_Bq \ar[rr]^-{\widetilde{\sigma}} \ar[d]_{\Delta_d\ubul_Bq} && q'\bul_Ac \ar[d]^{q'\ubul_Ac}\\
d\bul_Bd\bul_Bq \ar[r]_-{d\ubul_B\alpha} & d\bul_Bq\bul_Ac \ar[r]_-{\widetilde{\sigma}\ubul_Ac}
& q'\bul_Ac\bul_Ac
}\]
\[
\xymatrix{
d\bul_Bq \ar[rr]^-{\widetilde{\sigma}} \ar[d]_{\Delta_d\ubul_Bq} && q'\bul_Ac \ar[d]^{q'\ubul_Ac}\\
d\bul_Bd\bul_Bq \ar[r]_-{d\ubul_B\widetilde{\sigma}} & d\bul_Bq\bul_Ac \ar[r]_-{\alpha'\ubul_Ac}
& q'\bul_Ac\bul_Ac
}
\]
\item a $\dd\hbox{-}\cc$--bicomodule morphism
$\hat{\sigma}:d\bul_Bq\to d\bul_Bq'$, where the $\dd\hbox{-}\cc$--bicomodule structures on $d\bul_Bq$ and $d\bul_Bq'$ are obtained
from \leref{bicomodule}.
\end{enumerate}
\end{lemma}

\begin{proof}
Let us just give the corresponding formulae. The remaining part of the proof is just a computation. The equivalence of $(i)$ and $(ii)$ can also be deduced as a dualization of \cite[page 249]{LacStr}, and the equivalence of $(i)$ and $(iii)$ as a formalization of \cite[Proposition 2.2]{BrzElkGT}. Take $\sigma$ as in statement $(i)$, then define
$$\widetilde{\sigma}=(\sigma\ubul_Ac)\circ(d\ubul_B\alpha)\circ(\Delta_d\ubul_Bq),\qquad
\hat{\sigma}=(d\ubul_B\sigma)\circ(\Delta_d\ubul_Bq).$$
Conversely, if $\hat{\sigma}$ or $\widetilde{\sigma}$ are given, we can define
$$\sigma=(\varepsilon_d\ubul_Bq') \circ \widehat{\sigma},\qquad \sigma=(q'\ubul_A\varepsilon_c)\circ \widetilde{\sigma}\vspace{-17pt}.$$
\end{proof}

\begin{Definition}
We denote by $\REM(\Bb)$ the $2$--category whose $0$--cells and
$1$--cells are precisely those of $\RCOM(\Bb)$ and whose $2$--cells
are the $2$--cells $\sigma$ in $\Bb$ that satisfy the condition of
\leref{sigma} (i). Similarly one introduces the left-handed
version $\LEM(\Bb)$.
\end{Definition}

There exist a locally faithful pseudo functor
\begin{eqnarray}\eqlabel{forget2}
U:
&\REM(\Bb)&\to \Bb;\\
& \cc=(A,c,\Delta,\varepsilon)&\mapsto A\nonumber\\
& (q,\alpha) & \mapsto d\bul_Bq\nonumber\\
& \sigma & \mapsto \hat{\sigma}=(d\ubul_B\sigma)\circ(\Delta_d\ubul_Bq)\nonumber
\end{eqnarray}
Of course we can as well introduce the left hand versions of this functor,
$U_\ell: \LEM(\Bb)^{op}\to \Bb$.

\begin{theorem}\thlabel{localadjoint}
The pseudo functor $U$ from \equref{forget2} has locally a right adjoint.
\end{theorem}

\begin{proof}
Consider two comonads $\cc=(A,c,\Delta_c,\varepsilon_c)$ and $\dd=(B,d,\Delta_d,\varepsilon_d)$ in $\Bb$. The pseudo functor $U$ evaluated at $\dd$ and $\cc$ induces a functor
\[
\begin{array}{rrcl}
U_{\dd,\cc}:& \Hom_{\REM(\Bb)}(\dd,\cc) &\to& \Hom_\Bb(B,A)\\
& (q,\alpha) &\mapsto & d\bul_B q\\
\end{array}
\]
which has a right adjoint defined as follows
\[
\begin{array}{rrcl}
V'_{\dd,\cc}:&  \Hom_\Bb(B,A) &\to&\Hom_{\REM(\Bb)}(\dd,\cc)\\
& m & \mapsto & (m\bul_Ac,\varepsilon_d\ubul_Bm\ubul_B\Delta_c)
\end{array}
\]
The unit $\theta$ and counit $\epsilon$ of the adjunction are given as follows
\begin{eqnarray*}
\theta_q = (d\ubul_B\alpha)\circ(\Delta_d\ubul_Bq): d\bul_Bq \to d\bul_Bq\bul_Bc && \textrm{for all } (q,\alpha)\in\Hom_{\REM(\Bb)}(\dd,\cc)\\
\epsilon_m = m\ubul_A\varepsilon_c : m\bul_Ac\to m && \textrm{for all } m\in \Rcom(B,A)
\end{eqnarray*}
In view of \leref{bicomodule}, $\theta_q$ is given exactly by the right $\cc$--comodule structure on $d\bul_Bq$.
\end{proof}

Considering a trivial comonad $(\Omega,\Omega,\Omega,\Omega)$, as in \exref{trivcomonad}, we obtain from \thref{localadjoint} immediately the following well-known result.

\begin{corollary}\colabel{forgetinduction}
Let $\cc=(A,c,\Delta_c,\varepsilon_c)$ be a comonad in $\Bb$ and $\Omega$ any $0$--object in $\Bb$. Then the forgetful functor $\Rcom(\Omega,\cc)\to \Hom(\Omega,A)$ has a right adjoint given by $-\bul_Ac:\Hom(\Omega,A)\to \Rcom(\Omega,\cc)$.
\end{corollary}

Although the definition of a comonad is perfectly left-right
symmetric, we have two different possibilities for the definition
of the $2$--category of comonads, both in the original ($\LCOM$
and $\RCOM$) and in the modified case ($\LEM$ and $\REM$). This is
due to the asymmetry in the definition of comonad-morphisms. The
following proposition shows that there exists a `local' duality
between the left and right versions, however in general no duality
can be obtained for the whole $2$--categories.

Recall that an \emph{adjoint pair} in $\Bb$ is a sextuple
$\ppp=(A,B,p,q,\mu,\eta)$ where $A$ and $B$ are $0$--cells,
$p\in\Hom_1(A,B)$, $q\in\Hom_1(B,A)$,
$\mu\in{^A\Hom_2^A}(p\bul_Bq,A)$ and
$\eta\in{^B\Hom_2^B}(B,q\bul_Ap)$, such that the following
diagrams commute
\begin{equation}\eqlabel{ccc}
\xymatrix{
q \ar[rr]^-{\cong} \ar[d]_{\cong} && B\bul_Bq \ar[d]^{\eta\ubul_Bq}\\
q\bul_AA && q\bul_Ap\bul_Bq \ar[ll]^{q\ubul_A\mu} }\qquad
\xymatrix{
p \ar[rr]^-{\cong} \ar[d]_{\cong} && p\bul_BB \ar[d]^{p\ubul_B\eta}\\
A\bul_Ap && p\bul_Bq\bul_Ap \ar[ll]^{\mu\ubul_Aq} }
\end{equation}

If $(A,B,p,q,\mu,\eta)$ is an adjoint pair in $\Bb$, then we have
the following monad and comonad
$$(B,q\bul_Ap,q\ubul_A\mu\ubul_Ap,\eta),\qquad
(A,p\bul_Bq,q\ubul_B\eta\ubul_Bp,\mu)$$ More generally, given a
monad $\rrr=(A,r,\mu_r,\eta_r)$ and a comonad
$\dd=(B,d,\Delta_d,\varepsilon_d)$, we can construct from the
adjoint pair a new monad
\begin{equation}\eqlabel{adjointmonad}
(B,q\bul_Ar\bul_Ap,q\ubul_A\mu_r\ubul_Ap \circ
(q\ubul_Ar\ubul_A\mu\ubul_Ar\ubul_Ap),(q\ubul_A\eta_r\ubul_Ap)\circ
\eta)
\end{equation}
and comonad
\begin{equation}\eqlabel{adjointcomonad}
(A,p\bul_Bd\bul_Bq,(p\ubul_Bd\ubul_B\eta\ubul_Bd\ubul_Bq)\circ
p\ubul_B\Delta_d\ubul_Bq,\mu\circ(p\ubul_B\varepsilon_d\ubul_Bq)).
\end{equation}

The following proposition generalizes \cite[Proposition 2.1]{GT:comonad}. The equivalence between $(ii)$ and $(iii)$ is already stated in \cite[Theorem 9]{Str}.

\begin{proposition}\prlabel{locdual}
Suppose $\ppp=(A,B,p,q,\mu,\eta)$ is an adjoint pair, $\cc=(A,c,\Delta_c,\varepsilon_c)$ and $\dd=(B,d,\Delta_d,\varepsilon_d)$ are two comonads in $\Bb$. Then the following statements are equivalent
\begin{enumerate}[(i)]
\item There exists a homomorphism of $A$--comonads
$\varphi\in{^A\Hom_2^A}(p\bul_Bd\bul_Bq, c)$; \item there exists
an $\alpha\in{^B\Hom_2^A}(d\bul_Bq,q\bul_Ac)$, such that
$(q,\alpha)$ is a right comonad-morphism from $\dd$ to $\cc$;
\item there exists a $\beta\in{^A\Hom_2^B}(p\bul_Bd,c\bul_Ap)$,
such that $(p,\beta)$ is a left comonad-morphism from $\dd$ to
$\cc$.
\end{enumerate}
When these equivalent conditions hold, we say that $(p,q)$ is a
comonad-morphism with adjunction from $\dd$ to $\cc$.
\end{proposition}

\begin{proof}
$\ul{(i)\Rightarrow(ii)}$ Suppose there exists a morphism of comonads $\varphi$ as in the statement, then we define $\alpha:q\bul_Bd\to c\bul_Aq$ as follows
\[
\alpha:\xymatrix{d\bul_Bq\ar[rr]^-{\eta\ubul_Bd\ubul_Bq} && q\bul_Ap\bul_Bd\bul_Bq \ar[r]^-{q\ubul_A\varphi}
& q\bul_Ac.}
\]
Consider the following diagram.
\[
\xymatrix{
d\bul_Bq \ar[rr]^{\Delta_d\ubul_Bq} \ar[dd]_{\eta\ubul_Bd\ubul_Bq} && d\bul_Bd\bul_Bq \ar[d]^{d\ubul_B\eta\ubul_Bd\ubul_Bq}\\
&& d\bul_Bq\bul_Ap\bul_Bd\bul_Bq\ar[d]^{d\ubul_Bq\ubul_A\varphi}\\
q\bul_Ap\bul_Bd\bul_Bq \ar[r]^{q\ubul_Ap\ubul_B\Delta_d\ubul_Bq} \ar[dd]_{q\ubul_A\varphi}
& q\bul_Ap\bul_Bd\bul_Bd\bul_Bq \ar[d]^{q\ubul_A\ubul_Bd\ubul_B\eta\ubul_Bd\ubul_Bq}
& d\bul_Bq\bul_Ac \ar[d]^{\eta\ubul_Bd\ubul_Bq\ubul_Ac}\\
& q\bul_Ap\bul_Bd\bul_Bq\bul_Ap\bul_Bd\bul_Bq \ar[dr]_{q\ubul_A\varphi\ubul_A\varphi}
& q\bul_Ap\bul_Bd\bul_Bq\bul_Ac \ar[d]^{q\ubul_A\varphi\ubul_Ac}\\
q\bul_Ac \ar[rr]_{q\ubul_A\Delta_c} && q\bul_Ac\bul_Ac
}
\]
The outer diagram expresses the first condition for $\alpha$ to be
a right comonad-morphism. The upper part of the diagram commutes
by application of \equref{interchange2}, the lower part expresses
that $\varphi$ is a morphism of comonads. The second condition can
be computed as follows
\begin{eqnarray*}
(q\ubul_A\varepsilon_c)\circ \alpha&=&(q\ubul_A\varepsilon_c)\circ(q\ubul_A\varphi)\circ(\eta\ubul_Bd\ubul_Bq)\\
&=&(q\ubul_A\mu)\circ(q\ubul_Ap\ubul_B\varepsilon_d\ubul_Bq)\circ(\eta\ubul_Bd\ubul_Bq)\\
&=&(q\ubul_A\mu)\circ(\eta\ubul_Bq)\circ(\varepsilon_d\ubul_Bq)\\
&=&\varepsilon_d\ubul_Bq.
\end{eqnarray*}
Here we used in the second equality the counit condition of the morphism of comonads $\varphi$, third equation is again an application of \equref{interchange2} and the last equation follows from the fact that $\ppp$ is an adjoint pair in $\Bb$.

$\ul{(ii)\Rightarrow(i)}$ Suppose $\alpha$ exists as in the statement, then we define
\[
\varphi:\xymatrix{p\bul_Bd\bul_Aq \ar[r]^-{p\ubul_B\alpha} & p\bul_Bq\bul_Ac \ar[r]^-{\mu\ubul_Ac} & c.}
\]
We have to check that $\varphi$ is a morphism of comonads.
The following quadrangle expresses the counit condition on $\varphi$.
\[
\xymatrix{
p\bul_Bd\bul_Bq \ar[r]^{p\ubul_B\alpha} \ar[d]_{p\ubul_B\varepsilon_d\ubul_Bq}
& p\bul_Bq\bul_Ac \ar[r]^-{\mu\ubul_Ac} \ar[dl]^{p\ubul_Bq\ubul_A\varepsilon_c}
& c \ar[d]^{\varepsilon_c}\\
p\bul_Bq \ar[rr]_\mu && A
}
\]
The inner triangle commutes because of the counit condition on
$\alpha$, while the inner quadrangle commutes by
\equref{interchange2}. Next, we verify the compatibility of
$\varphi$ with the comultiplication maps.
\begin{eqnarray*}
(\varphi\ubul_A\varphi)\circ(\Delta_{p\bul_Bd\bul_Bq})
& =&(\mu\ubul_Ac\ubul_A\mu\ubul_Ac)\circ(p\ubul_B\alpha\ubul_B\alpha)\circ(p\ubul_Bd\ubul_B\eta\ubul_Bd\ubul_Bq)\circ(p\ubul_B\Delta_d\ubul_Bq)\\
&=&(\mu\ubul_Ac\ubul_Ac)\circ(p\ubul_Bq\ubul_Ac\ubul_A\mu\ubul_Ac)\circ(p\ubul_B\alpha\ubul_Ap\ubul_Bq\ubul_Ac)\\
&&\circ(p\ubul_Bd\ubul_Bq\ubul_Ap\ubul_B\alpha)\circ(p\ubul_Bd\ubul_B\eta\ubul_Bd\ubul_Bq)\circ(p\ubul_B\Delta_d\ubul_Bq)\\
&=&(\mu\ubul_Ac\ubul_Ac)\circ(p\ubul_B\alpha\ubul_Ac)\circ(p\ubul_Bd\ubul_Bq\ubul_A\mu\ubul_Ac)\\
&&\circ(p\ubul_Bd\ubul_B\eta\ubul_Bq\ubul_Ac)\circ(p\ubul_Bd\ubul_B\alpha)\circ(p\ubul_B\Delta_d\ubul_Bq)\\
&=&(\mu\ubul_Ac\ubul_Ac)\circ(p\ubul_B\alpha\ubul_Ac)\circ(p\ubul_Bd\ubul_B\alpha)\circ(p\ubul_B\Delta_d\ubul_Bq)\\
&=&(\mu\ubul_Ac\ubul_Ac)\circ(p\ubul_Bq\ubul_A\Delta_c)\circ(p\ubul_B\alpha)\\
&=&\Delta_c\circ(\mu\ubul_Ac)\circ(p\ubul_B\alpha)\\
&=&\Delta_c\circ\varphi
\end{eqnarray*}
The second, third and penultimate equality are applications of
\equref{interchange2}, the fourth one follows from the condition
on the adjoint pair $\ppp$, and the fifth equation uses the fact
that $\alpha$ is a right comonad-morphism.

The equivalence $\ul{(i)\Leftrightarrow(iii)}$ follows by left-right duality from $(i)\Leftrightarrow(ii)$.
\end{proof}

\begin{Definition}
The $2$--category $\fREM(\Bb)$ consists of the following objects
\begin{itemize}
\item  $0$--cells are comonads; \item a $1$--cell from $\cc$ (over
$A$) to $\dd$ (over $B$) consists of an adjoint pair
$(A,B,p,q,\mu,\eta)$ in $\Bb$ together with a right
comonad-morphism from $\dd$ to $\cc$ of the form $(q,\alpha)$;
\item a $2$--cell in $\fREM(\Bb)$ is a $2$--cell $\sigma$ in $\Bb$
that satisfies the conditions of \leref{sigma} (i). \item To
define the composition of $1$--cells, we use the composition of
comonad-morphisms as in $\RCOM(\Bb)$ together with the usual
composition of adjoint pairs.
\end{itemize}
And similar for $\fLEM(\Bb)$.
\end{Definition}

\begin{remark}
Consider two $1$--cells in $\fREM(\Bb)$ from $\dd$ to $\cc$,
consisting of the same comonad-morphism $(q,\alpha)$ together with
an adjoint pair $(A,B,p,q,\mu,\eta)$, respectively
$(A,B,p',q,\mu',\eta')$. Since every adjoint pair
$(A,B,p,q,\mu,\eta)$ induces a pair of adjoint functors $-\bul_Bq
\dashv -\bul_Ap$ between the categories $\Hom_\Bb(\Omega,A)$ and
$\Hom_\Bb(\Omega,B)$ for each $0$--cell $\Omega$ in $\Bb$, we find
by the uniqueness of adjoint functors a natural isomorphism
$\theta:-\bul_Ap\to-\bul_Ap'$ of functors from
$\Hom_\Bb(\Omega,A)$ to $\Hom_\Bb(\Omega,B)$. Taking $\Omega=A$,
we find that $\theta_A:p\to p'$ is an isomorphism. Thus, the
adjoint $1$--cell of a comonad-morphism $(q,\alpha)$ is, up to
isomorphism, uniquely determined in $\Hom_\Bb(A,B)$.~\qed
\end{remark}

From \prref{locdual} we obtain the following generalization of
\cite[Proposition 2.5]{BrzElkGT} to our setting.

\begin{proposition}
There is a biequivalence between $\fREM(\Bb)$ and $\fLEM(\Bb)^{co}$.
\end{proposition}

\begin{proof}
Consider a $1$--cell from $\dd$ to $\cc$ in $\fREM(\Bb)$ consisting of a right comonad-morphism $(q,\alpha)$, together with an adjoint pair $(A,B,p,q,\mu,\eta)$. Then we know by \prref{locdual} that we can construct a left comonad-morphism $(p,\beta)$, hence a $1$--cell in $\fLEM(\Bb)$.

If we take now a $2$--cell $\sigma$ in $\fREM(\Bb)$ connecting the
$1$--cell $(q,\alpha)$ with $(A,B,p,q,\mu,\eta)$ and $(q',\alpha')$
with $(A,B,p',q',\mu',\eta')$, then $\sigma$ is a $2$--cell in
$\Bb$ of the form $d\ubul_B q\to q'$. We define now $\tau$ as a
$2$--cell in $\fLEM(\Bb)$ from the left comonad-morphism
$(p,\beta)$ to $(p',\beta')$ as follows:
\[
\xymatrix{
\tau:p'\bul_B d \cong p'\bul_Bd\bul_BB \ar[rr]^-{p\ubul_Bd\ubul_B\eta} && p'\bul_Bd\bul_Bq\bul_Ap
\ar[r]^-{p'\ubul_B\sigma\ubul_Ap} & p'\bul_Bq'\bul_Ap \ar[r]^-{\mu'\ubul_Ap} & A\bul_Ap\cong p.
}
\]
The remaining details are left to the reader.
\end{proof}

Consider $1$--cells $m,n:A\to B$ and $2$--cells $\alpha,\beta: m\to
n$. Suppose that the equalizer $(e,\epsilon)$ of $\alpha$ and
$\beta$ exists in $\Hom(A,B)$,
\[
\xymatrix{ e \ar[r]^\epsilon & m \ar@<.5ex>[r]^\alpha
\ar@<-.5ex>[r]_\beta & n. }
\]
Let $a:B\to C$ be any $1$--cell. We say that the equalizer
$(e,\epsilon)$ is right $a$--pure if and only if the following
diagram is an equalizer in $\Hom(A,C)$:
\[
\xymatrix{
e\bul_B a \ar[r]^{\epsilon\ubul_B a} & m\bul_B a \ar@<.5ex>[r]^{\alpha\ubul_Ba} \ar@<-.5ex>[r]_{\beta\ubul_Ba} & n\bul_Ba.
}
\]
If an equalizer is right $a$--pure for all choices of $a$, then we
just say that the equalizer is right pure. Similarly, let $b:C\to
A$ be a $1$--cell. We say that the equalizer $(e,\epsilon)$ is
left $b$--pure if and only if the following diagram
\[
\xymatrix{
b\bul_A e \ar[r]^{b\ubul_A \epsilon} & b\bul_Am \ar@<.5ex>[r]^{b\ubul_A\alpha} \ar@<-.5ex>[r]_{b\ubul_A\beta} & b\bul_An
}
\]
is an equalizer in $\Hom(C,A)$.

\begin{lemma}\lelabel{equalizer}
Consider a comonad $\dd=(B,d,\Delta_d,\varepsilon_d)$ and a
$0$--cell $\Omega$ in $\Bb$. Take $f,g : (m,\rho^m)\to
(n,\rho^n)\in \Rcom(\Omega,\dd)$. Suppose that the equalizer
$(q,\epsilon)$ of the pair $(f,g)$ exists in $\Hom(\Omega,B)$
(i.e. after applying the forgetful functor
$U_{\Omega,\dd}:\Rcom(\Omega,\dd)\to\Hom(\Omega,B)$). Then the
following statements are equivalent
\begin{enumerate}[(i)]
\item $(q,\epsilon)$ is right $d\bul_Bd$--pure;
\item the following equalizers exist in $\Rcom(\Omega,\dd)$ and are preserved by the forgetful functor $U_{\Omega,\dd}:\Rcom(\Omega,\dd)\to\Hom(\Omega,B)$.
\begin{eqnarray}
&\xymatrix{
q \ar[rr]^{\epsilon}
&& m \ar@<.5ex>[rr]^-{f} \ar@<-.5ex>[rr]_-{g}
&& n
}\\
&\xymatrix{
q\bul_Bd \ar[rr]^{\epsilon\ubul_Bd}
&& m\bul_Bd \ar@<.5ex>[rr]^-{f\ubul_Bd} \ar@<-.5ex>[rr]_-{g\ubul_Bd}
&& n\bul_Bd
}\eqlabel{equald} \\
&\xymatrix{
q\bul_Bd\bul_Bd \ar[rr]^{\epsilon\ubul_Bd\bul_Bd}
&& m\bul_Bd\bul_Bd \ar@<.5ex>[rr]^-{f\ubul_Bd\bul_Bd} \ar@<-.5ex>[rr]_-{g\ubul_Bd\bul_Bd}
&& n\bul_Bd\bul_Bd
} \eqlabel{dbuld}
\end{eqnarray}
\end{enumerate}
\end{lemma}

\begin{proof}
We only prove the implication $(i)\Rightarrow(ii)$, the converse
is easy. First remark that if an equalizer
$\xymatrix{q\ar[r]^\epsilon & m \ar[r]^{f,g} & n}$ is right
$d\bul_Bd$--pure, then this equalizer is also right $d$--pure.
This follows easily from the following diagram, once we observe
that $(d\bul_B\varepsilon_d) \circ \Delta_d=d$:
\[
\xymatrix{
q\bul_Bd \ar[rr]^{\epsilon\ubul_Bd} \ar@<.5ex>[d]^{q\ubul_B\Delta_d}
&& m\bul_Bd \ar@<.5ex>[rr]^-{f\ubul_Bd} \ar@<-.5ex>[rr]_-{g\ubul_Bd} \ar@<.5ex>[d]^{m\ubul_B\Delta_d}
&& n\bul_Bd \ar@<.5ex>[d]^{\Delta_d\ubul_Bn} \\
q\bul_Bd\bul_Bd \ar[rr]^{\epsilon\ubul_Bd\ubul_Bd} \ar@<.5ex>[u]^{q\ubul_B\varepsilon_d\ubul_Bd}
&& m\bul_Bd\bul_Bd \ar@<.5ex>[rr]^-{f\ubul_Bd\ubul_Bd} \ar@<-.5ex>[rr]_-{g\ubul_Bd\ubul_Bd} \ar@<.5ex>[u]^{m\ubul_B\varepsilon_d\ubul_Bd}
&& n\bul_Bd\bul_Bd \ar@<.5ex>[u]^{n\ubul_B\varepsilon_d\ubul_Bd}
}
\]

We know by assumption that the equalizer \equref{equald} exists
in $\Hom(\Omega,B)$. Moreover, the functor
$-\bul_Bd:\Hom(\Omega,B)\to\Rcom(\Omega,\dd)$ preserves
equalizers, since it is the right adjoint to the forgetful functor
$U_{\Omega,\dd}$. Therefore the equalizer \equref{equald} exists
in $\Rcom(\Omega,\dd)$ and is preserved by the forgetful functor.
Applying this argument a second time we obtain the same result for
the equalizer \equref{dbuld}.

From the $d$--purity we obtain the following commutative diagram, where horizontal rows are equalizers.
\[
\xymatrix{
q \ar[rr]^{\epsilon} \ar@{.>}[d]
&& m \ar@<.5ex>[rr]^-{f} \ar@<-.5ex>[rr]_-{g} \ar[d]^{\rho^m}
&& n \ar[d]^{\rho^n} \\
q\bul_Bd \ar[rr]^{\epsilon\ubul_Bd}
&& m\bul_Bd \ar@<.5ex>[rr]^-{f\ubul_Bd} \ar@<-.5ex>[rr]_-{g\ubul_Bd}
&& n\bul_Bd
}
\]
From the universal property of the equalizer, we therefore obtain
a $2$--cell $\rho:q\to q\bul_Bd$. One can easily prove that
$(q\bul_B\varepsilon_d)\circ\rho=q$ and the coassociativity of
this coaction follows form the right $d\bul_Bd$--purity.

Let us prove that $((q,\rho),\epsilon)$ is an equalizer in $\Rcom(\Omega,\dd)$. Suppose there exists a morphism $\kappa:(q',\rho')\to(m,\rho^m)\in\Rcom(\Omega,\dd)$ such that $f\circ\kappa=g\circ\kappa$. Then we know that there exists a $2$--cell $\lambda\in{^\Omega\Hom^B_2}(q',q)$ such that $\kappa=\epsilon\circ\lambda$. Let us check that $\lambda$ is right $\dd$--colinear. We use the defining property of $\rho$ in the first equation, the relation between $\kappa$ and $\lambda$ in the second and last equality and the right $\dd$--colinearity of $\kappa$ in the third one.
\begin{eqnarray*}
(\epsilon\ubul_Bd)\circ\rho\circ\lambda &=& \rho^m\circ\epsilon\circ\lambda=\rho^m\circ\kappa\\
&=&(\kappa\ubul_Bd)\circ\rho'=(\epsilon\ubul_Bd)\circ(\lambda\ubul_Bd)\circ\rho'.
\end{eqnarray*}
From this we obtain that $\rho\circ\lambda=(\lambda\ubul_Bd)\circ\rho'$ since $\epsilon\ubul_Bd$ is a monomorphism in $\Hom(\Omega,B)$, as $(q\bul_Bd,\epsilon\ubul_Bd)$ is an equalizer.
\end{proof}

\begin{theorem}
Consider a comonad $\dd=(B,d,\Delta_d,\varepsilon_d)$ and a $0$--cell $\Omega$ in $\Bb$.
If $\Hom(\Omega,B)$ has all equalizers, then the following statements are equivalent.
\begin{enumerate}[(i)]
\item $T = - \bul_B d : \Hom (\Omega,B) \rightarrow
\Hom(\Omega,B)$ preserves equalizers; \item $\Rcom(\Omega,\dd)$
has all equalizers, and they are preserved by the forgetful
functor $U_{\Omega,\dd}:\Rcom(\Omega,\dd)\to\Hom(\Omega,B)$; \item
the forgetful functor
$U_{\Omega,\dd}:\Rcom(\Omega,\dd)\to\Hom(\Omega,B)$ preserves
equalizers.
\end{enumerate}
\end{theorem}

\begin{proof}
$\ul{(i)\Rightarrow(ii)}$. Consider two parallel morphisms $f,g:(m,\rho^m)\to(n,\rho^n)$ in $\Rcom(\Omega,\dd)$, and construct the equalizer $(q,\epsilon)$ of $(f,g)$ in $\Hom(\Omega,\dd)$. Since $T$ preserves all equalizers, this equalizer is $d\bul_Bd$--pure. Therefore, we can apply \leref{equalizer} and we obtain that $(q,\epsilon)$ is also an equalizer in $\Rcom(\Omega,\dd)$.\\
$\ul{(ii)\Rightarrow(iii)}$. Trivial.\\
$\ul{(iii)\Rightarrow(i)}$. Consider $2$--cells
$f,g\in{^\Omega\Hom^B_2}(m,n)$, and let $(q,\epsilon)$ be the
equalizer of $(f,g)$ in $\Hom(\Omega,B)$. The functor
$-\bul_Bd:\Hom(\Omega,B)\to\Rcom(\Omega,\dd)$ preserves
equalizers, since it is the right adjoint to the forgetful functor
$U_{\Omega,\dd}$. As the forgetful functor preserves equalizers as
well, we find that $T = U_{\Omega,\dd}\circ (-\bul_Bd)$ preserves
equalizers.
\end{proof}

Consider now a comonad $\cc=(A,c,\Delta_c,\varepsilon_c)$, a right
$\cc$--comodule of $C$--type $(m,\rho^m)$ and a left
$\cc$--comodule of $B$--type $(n,\lambda^n)$ for all $0$--cells
$B$ and $C$. If it exists, then we will denote the equalizer of
$\rho^m\ubul_An$ and $m\ubul_A\lambda^n$ (in $\Hom(C,B)$) by
$m\bul^c n$ and we call this the \emph{cotensor}\index{cotensor
product} product,
\begin{equation}\eqlabel{cotensorprod}
\xymatrix{
m\bul^c n \ar[r] & m\bul n \ar@<.5ex>[r]^-{m\ubul_A\lambda^n} \ar@<-.5ex>[r]_-{\rho^m\ubul_An} & m\bul_Ac\bul_An.
}
\end{equation}
Consider now two other comonads $\dd=(B,d,\Delta_d,\varepsilon_d)$ and $\ee=(C,e,\Delta_e,\varepsilon_e)$.
Supose that $(n,\lambda^n,\rho^n)$ is a $\cc$-$\dd$--bicomodule and $(m,\lambda^m,\rho^m)$ is an $\ee$-$\cc$--bicomodule.

\begin{corollary}\colabel{cotensorcomodule}
With notation as introduced above, suppose that the equalizer $m\bul^cn$ exists in $Hom(C,B)$.
\begin{enumerate}[(i)]
\item If the cotensor product $m\bul^c n$ is right
$d\bul_Bd$--pure, then $m\bul^c n$ is a right $\dd$--comodule. \item
If the cotensor product $m\bul^c n$ is left $e\bul_Ce$--pure, then
$m\bul^c n$ is a left $\ee$--comodule. \item If the cotensor
product $m\bul^c n$ is both right $d\bul_Bd$--pure and left
$e\bul_Ce$--pure, then $m\bul^c n$ is a  $\ee$-$\dd$--bicomodule.
\end{enumerate}
\end{corollary}

Consider a $2$--category $\Bb$ such that for every pentuple of
comonads $\cc$ (over $A$), $\dd$ (over $B$), $\ee$ (over $C$),
$\ff$ (over $D$) and $\gg$ (over $E$) and for all
$\ee$-$\cc$--bicomodule $m$ and every $\cc$-$\dd$--bicomodule $n$,
the cotensor product $m\bul^c n$ exists and is left $p$--pure and
right $q$--pure for any $\ff$-$\ee$--bicomodule $p$ and any
$\dd$-$\gg$--bicomodule $q$. We can now construct a new bicategory
$\Bic(\Bb)$ out of $\Bb$:
\begin{itemize}
\item A $0$--cell in $\Bic(\Bb)$ is a comonad in $\Bb$;
\item a
$1$--cell $\cc\to\dd$ in $\Bic(\Bb)$ is a $\cc$-$\dd$--bicomodule;
\item a $2$--cell $\varphi$ in $\Bic(\Bb)$ between two
$\cc$-$\dd$--bicomodules, is a left $\cc$-- and a right
$\dd$--colinear morphism.
\end{itemize}
Composition of $1$--cells is given by the cotensor product, unit elements are the comonads, considered as bicomodules over themselves.

\begin{remark}
It follows from \coref{cotensorcomodule} that the composition of $1$--cells in $\Bic(\Bb)$ is well-defined. Moreover the composition of $1$-cells is also associative (up to isomorphism), as we can see as follows. Consider four comonads in $\Bb$: $\cc$ (over $A$), $\dd$ (over $B$), $\ee$ (over $C$), $\ff$ (over $D$) together with a $\cc$-$\dd$--bicomodule $m$, a $\dd$-$\ee$--bicomodule $n$ and a $\ee$-$\ff$--bicomodule $p$. Then we need an isomorphism of the form
$$(m\bul^d n)\bul^e p\cong m\bul^d (n\bul^e p).$$
Consider the following diagram
\[
\xymatrix{
(m\bul^d n)\bul^e p \ar[r] \ar[d]^{\psi_1} & (m\bul_B n)\bul^e p \ar@<.5ex>[r] \ar[d]^-{\psi_2} \ar@<-.5ex>[r] & (m\bul_Bd\bul_B n)\bul^e p \ar[d]^-{\psi_3}\\
m\bul^d (n\bul^e p) \ar[r] & m\bul_B (n\bul^e p) \ar@<.5ex>[r]^-{} \ar@<-.5ex>[r]^-{}
& m\bul_Bd\bul_B (n\bul^e p)
}
\]
The bottom row is exact by the definition of the cotensor product. The exactness of the top row can be deduced from the fact that $m\bul^d n$ is right $p$--, right $e$-- and right $e\bul_Cp$--pure. Furthermore we know by the left $m$-- and left $m\bul_B d$--purity of $n\bul^e p$ that $\psi_2$ and $\psi_3$ are isomorphisms, hence $\psi_1$ is an isomorphism as well by the univeral property of the equalizer.

An example of a bicategory $\Bb$ as above is given by the bicategory of bimodules over division algebras.
\end{remark}

\begin{theorem}\thlabel{REMBic}
Let $\Bb$ be a $2$--category as above. Then there exist locally
full and locally faithful pseudo functors
$$F:\REM(\Bb)\to \Bic(\Bb), \qquad G:\LEM(\Bb)^{co}\to \Bic(\Bb).$$
\end{theorem}

\begin{proof}
Since both the $0$--cells in $\REM(\Bb)$ and $\Bic(\Bb)$ are
comonads, we can define $F$ to be the identity on $0$--cells. A
$1$--cell in $\REM(\Bb)$ is a right comonad-morphism
$(q,\alpha):\dd\to \cc$. By \leref{bicomodule}, $d\bul_Bq$ is a
$\dd$-$\cc$--bicomodule, thus we define $F(q,\alpha)=d\bul_Bq$. A
$2$--cell in $\REM(\Bb)$ is a comonad transformation
$\sigma:(q,\alpha)\to(q',\alpha')$. We know from \leref{sigma}
that the second unreduced form $\hat{\sigma}$ of $\sigma$ is a
bicomodule morphism from $F(q,\alpha)=d\bul_Bq$ to
$F(q',\alpha')=d\bul_Bq'$, thus we define
$F(\sigma)=\hat{\sigma}$.

Finally we have to find natural isomorphisms $\gamma$ \equref{laxfunctor} and $\delta$ \equref{laxfunctor2}. To this end, consider a comonad-morphism $(q,\alpha): \ee\to \dd$  and a comonad-morphism $(q',\alpha'): \dd\to\cc$. We need to prove that $F(q)\bul^dF(q')\cong F(q\bul_Bq')$. Indeed:
$$(e\bul_C q)\bul^d(d\bul_B q') \cong ((e\bul_C q)\bul^dd) \bul_B q'
\cong e\bul_Cq\bul_Bq'$$
where the first isomorphism is a consequence from the purity conditions.
The property for $\delta$ is trivial, the unit objects in both categories are comonads.

That the pseudo functors $F$ and $G$ are locally full and locally faithful follows from the description of $2$--cells in $\REM(\Bb)$ as bicomodule-morphisms in \leref{sigma}.
\end{proof}

\section{Equivalences between Comodule categories}\selabel{Galois}

In this section we will show that the theory developed in
\cite{BrzElkGT} for the bicategory  $\Bim$ of bimodules, can be
established in any $2$--category $\Bb$, and, since every
bicategory is biequivalent to a suitable $2$--category, it is
already extended to any bicategory (see \seref{definitions}). In
fact, the results of \cite[\S 5]{BrzElkGT} are improved in what
follows even in the case $\Bb = \Bim$.

\subsection{Push-out and pull-back functors}\selabel{popb}

Given comonads $\dd=(B,d,\Delta_d,\varepsilon_d)$ and
$\cc=(A,c,\Delta_c,\varepsilon_c)$ in the $2$--category $\Bb$,
every comonad-morphism $(q,\alpha)$ from $\dd$ to $\cc$ defines,
for each object $\Omega$ of $\Bb$, a functor
$$\qQ:\Rcom(\Omega,\dd)\to\Rcom(\Omega,\cc),$$
given by
$\qQ(m,\rho)=(m\bul_Bq,(m\ubul_B\alpha)\circ(\rho\ubul_Bq))$, for
any comodule $(m,\rho)\in\Rcom(\Omega,\dd)$ and $\qQ(\phi)=\phi\ubul_Bq$
for any morphism of comodules $\phi\in\Rcom(\Omega,\dd)$. This functor
makes commute the diagram
\begin{equation}\eqlabel{forgetdiagram} \xymatrix{
\Rcom(\Omega,\dd)\ar[d] \ar[d]_{U_{\Omega,\dd}} \ar[r]^-\qQ &\Rcom(\Omega,\cc) \ar[d]^{U_{\Omega,\cc}} \\
\Hom(\Omega,B)\ar[r]_-{- \bul_B q} & \Hom(\Omega,A).
}\end{equation}

Conversely, if $\qQ : \Rcom(B,\dd) \rightarrow \Rcom(B,\cc)$ is a
functor making commute \equref{forgetdiagram} for $\Omega = B$ and
some $1$--cell $q : B \rightarrow A$, then there is a
comonad-morphism $(q,\alpha)$ which induces $\qQ$. In fact,
$\qQ(d) = d \bul_B q$ as an object in $\Hom(B,A)$.  Therefore, we
find that $d\bul_Bq$ possesses a right $\cc$--comodule structure
and, from the functoriality of $\qQ$, it follows that
$\qQ(\Delta_d)$ must be a homomorphism of right $\cc$--comodules,
hence $d\bul_Bq$ is even a $\dd$-$\cc$--bicomodule. By
\leref{bicomodule} we obtain a comonad-morphism from $\dd$ to
$\cc$ of the form $(q,\alpha)$.

\begin{Definition}
Inspired by the terminology of \cite{BrzElkGT}, we will call the
functor $\qQ$ the \emph{pushout} functor associated to the comonad
morphism $(q,\alpha)$.
\end{Definition}
It is a natural question to pose whether the pushout functor has a
right adjoint. Generalizing results from \cite{BrzElkGT} and
\cite{GT:comonad}, we find a criterion for this to hold in the
`finite case'.
Consider a left comonad-morphism $(p,\beta)$ from $\dd$ to $\cc$.
We know from \leref{bicomodule} that $p\bul_Bd$ is a
$\cc\hbox{-}\dd$--bicomodule. In this note we will say that the
category $\Rcom(\Omega,\dd)$ satisfies the \emph{equalizer
condition for $p$} if for all comodules $(n,\rho)$ in
$\Rcom(\Omega,\cc)$ the equalizer of $\rho\ubul_Ap\ubul_Bd$ and
$n\ubul_A\lambda^{p\bul_Bd}$ exists in $\Rcom(\Omega,\dd)$. In
view of \coref{cotensorcomodule}, a sufficient condition for this
to be satisfied is that the following cotensor product exists in
$\Hom(\Omega,B)$ and is right $d\bul_Bd$--pure.
\begin{equation}\eqlabel{cotensor}
\xymatrix{n\bul^c(p\bul_Bd)\ar[r]^{{\rm eq}_n} &
n\bul_A(p\bul_Bd)\ar@<-.5ex>[rr]_-{n\ubul_A\lambda^{p\bul_Bd}}
\ar@<.5ex>[rr]^-{\rho\ubul_Ap\ubul_Bd} &&
n\bul_Ac\bul_A(p\bul_Bd).}
\end{equation}
For this reason, we will denote the equalizer of
$\rho\ubul_Ap\ubul_Bd$ and $n\ubul_A\lambda^{p\bul_Bd}$ in
$\Rcom(\Omega,\dd)$ by
$(n\bul^c(p\bul_Bd),\rho^{n\bul^c(p\bul_Bd)})$, if
$\Rcom(\Omega,\dd)$ satisfies the equalizer condition.

\begin{proposition}
Consider comonads $\dd$ and $\cc$ in $\Bb$ and take
$p\in\Hom_1(A,B)$.
\begin{enumerate}
\item Let $\Omega$ be a $0$--cell in $\Bb$ and $(p,\beta)$ a left
comonad-morphism from $\dd$ to $\cc$. If $\Rcom(\Omega,\dd)$
satisfies the equalizer condition for $p$, then there exists a
functor
$$\pP:\Rcom(\Omega,\cc)\to\Rcom(\Omega,\dd),$$
such that the following diagram commutes
\begin{equation}\eqlabel{inductiondiagram}
\xymatrix{
\Rcom(\Omega,\cc)  \ar[r]^-\pP &\Rcom(\Omega,\dd) \\
\Hom(\Omega,A)\ar[r]_-{-\bul_Ap} \ar[u]^{-\bul_Ac} &
\Hom(\Omega,B) \ar[u]_{-\bul_Bd} ;}
\end{equation}
\item Conversely, if a functor $\pP$ as in part (1) exists for the
case $\Omega=A$, then there exists a left comonad-morphism from
$\dd$ to $\cc$ of the form $(p,\beta)$.
\end{enumerate}
\end{proposition}

\begin{proof}
\ul{(1)}. We define the functor $\pP$ by $\pP(n,\rho)=(n\bul^c(p\bul_Bd),\rho^{n\bul^c(p\bul_Bd)})$ for all comodules $(n,\rho)$ in $\Rcom(\Omega,\cc)$.\\
\ul{(2)}. Consider the element $A\in\Hom(A,A)$, by the
commutativity of diagram \equref{inductiondiagram} we can compute
$\pP(A\bul_Ac)\cong\pP(c)\cong p\bul_Bd$ as an object in
$\Hom(A,B)$. Moreover, we find that $p\bul_Bd$ has a right
$\dd$--comodule structure induced by $\dd$ and a left
$\cc$--comodule structure induced by $\cc$ and the functoriality of
$\pP$. It is even a $\cc$-$\dd$--bicomodule. Hence we find by
\leref{bicomodule} a comonad-morphism of the form $(p,\beta)$.
\end{proof}

\begin{Definition}
The functor $\pP$ is called the \emph{pullback} functor associated
to the left comonad-morphism $(p,\beta)$.
\end{Definition}

Consider the situation where $(p,q):\dd\to\cc$ is a
comonad-morphism with adjunction. By \prref{locdual}, the
existence of the comonad-morphism with adjunction $(p,q)$ implies
the existence of right and left comonad-morphisms $(q,\alpha)$ and
$(p,\beta)$. Then we obtain for any $0$--cell $\Omega$ in $\Bb$ the
following diagram of functors
\begin{equation}\eqlabel{diagramGal}
\xymatrix{
\Rcom(\Omega,\cc) \ar@<-.5ex>@{.>}[rr]_-{\pP} \ar@<-.5ex>[d]_{U_{\Omega,\cc}}
&& \Rcom(\Omega,\dd) \ar@<-.5ex>[ll]_-{\qQ} \ar@<.5ex>[d]^{U_{\Omega,\dd}}\\
\Hom(\Omega,A) \ar@<.5ex>[rr]^-{-\bul_Ap}
\ar@<-.5ex>[u]_{-\bul_Ac} && \Hom(\Omega,B)
\ar@<.5ex>[ll]^-{-\bul_Bq} \ar@<.5ex>[u]^{-\bul_Ad}, }
\end{equation}
where the lower and the vertical pairs of functors are adjoint.
The dotted arrow in the upper row only exists if the equalizer
condition is satisfied, in that situation the upper pair is also a
pair of adjoint functors, as we will show in the next proposition.
Moreover, the diagram commutes in the following sense: the outer
square and the inner square both commute.

\begin{proposition}\prlabel{adjuntos}
Consider a comonad-morphism with adjunction $(p,q):\dd\to\cc$. If
$\Rcom(\Omega,\dd)$ satisfies the equalizer condition for $p$,
then the pullback functor $\pP$ associated to $p$ is a right
adjoint to the pushout functor $\qQ$ associated to $q$.
\end{proposition}

\begin{proof}
The adjoint pair $(A,B,p,q,\mu,\eta)$ induces an adjoint pair of
functors
\[
\xymatrix{ \Hom(\Omega,B) \ar@<.5ex>[rr]^-{-\bul_Bq} &&
\Hom(\Omega,A) \ar@<.5ex>[ll]^-{-\bul_Ap} }
\]
Let us consider, for any $x\in\Hom_1(\Omega,B)$ and
$y\in\Hom_1(\Omega,A)$, the adjointness isomomorphism
$$\Phi^0_{x,y}:{^\Omega\Hom_2^A}(x\bul_Bq,y)\to {^\Omega\Hom_2^B}(x,y\bul_Ap).$$

On the other hand, $-\bul_Bd$ defines a comonad on the category
$\Hom(\Omega,B)$, so that we obtain an adjoint pair
\[
\xymatrix{ \Rcom(\Omega,\dd) \ar@<.5ex>[rr]^-{U} && \Hom(\Omega,B)
\ar@<.5ex>[ll]^-{-\bul_Bd} },\] where $U$ denotes the forgetful
functor. The adjunction $U \dashv -\bul_Bd$ entails a second
natural isomorphism
$$\Phi^1_{x,y}:{^\Omega\Hom_2^B}(x,z)\to {^\Omega\Hom^\dd}(x,z\bul_Bd),$$
where $x,z\in\Hom_1(\Omega,B)$.

Finally we can construct the following diagram.
\[
\xymatrix{
{^\Omega\Hom^\cc}(x\bul_Bq,y) \ar@{.>}[rrr] \ar[d] &&& {^\Omega\Hom^\dd}(x,y\bul^c(p\bul_Bd)) \ar[d]\\
{^\Omega\Hom_2^A}(x\bul_Bq,y) \ar@<-.5ex>[d]_{\rho^y\circ-} \ar@<.5ex>[d]^{-\circ\rho^{x\ubul_Bq}} \ar[r]^{\Phi^0_{x,y}} &
{^\Omega\Hom_2^B}(x,y\bul_Ap) \ar[rr]^-{\Phi^1_{x,y\bul_Ap}} &&
{^\Omega\Hom^\dd}(x,y\bul_Ap\bul_Bd) \ar@<-.5ex>[d]_{(\rho^y\ubul_Ap\ubul_Bd)\circ-} \ar@<.5ex>[d]^{(y\ubul_A\lambda^{p\bul_Bd})\circ-} \\
{^\Omega\Hom_2^A}(x\bul_Bq,y\bul_Ac) \ar[r]_-{\Phi^0_{x,y\bul_Ac}}
& {^\Omega\Hom_2^B}(x,y\bul_Ac\bul_Ap)
\ar[rr]_-{\Phi^1_{x,y\bul_Ac\bul_Ap}}
&&{^\Omega\Hom^\dd}(x,y\bul_Ac\bul_Ap\bul_Bd) }
\]
Now check that the vertical lines are equalizers and that the
diagram commutes (following the equally aligned vertical lines).
Since both lower horizontal arrows are isomorphisms, we find the
existence of an isomorphism in the upper horizontal line as well,
by the universal property of the equalizers. This implies the
adjunction between the pushout functor $-\bul_Bq$ and the pullback
functor $-\bul^c(p\bul_Bd)$.
\end{proof}

We state the explicit form of the unit $\zeta$ and counit $\nu$ of
the adjunction $\qQ \dashv \pP$. To obtain $\nu$, consider the
equalizer from \equref{cotensor},
$$
\xymatrix{y\bul^c(p\bul_Bd)\ar[r]^{{\rm eq}_y} & y\bul_A(p\bul_Bd)\ar@<-.5ex>[r] \ar@<.5ex>[r] & y\bul_Ac\bul_A(p\bul_Bd)}
$$
and apply the pushout functor $-\bul_Bq$ on this exact row, then we obtain
\begin{equation}\eqlabel{nu}
\xymatrix{(y\bul^c(p\bul_Bd))\ar@{.>}[rd]_{\nu_y} \bul_Bq \ar[r]^{{\rm eq}_y\ubul_Bq}
 & y\bul_A(p\bul_Bd)\bul_Bq\ar@<-.5ex>[r]\ar[d]^{y\ubul_A\epsilon} \ar@<.5ex>[r] &y\bul_Ac\bul_A(p\bul_Bd)\bul_Bq.\\
& y}
\end{equation}
So $\nu$ is given by the formula
$\nu_y=(y\bul_A\epsilon)\circ({\rm eq}\bul_Bq)$, where $\epsilon$
denotes the counit of the comonad $p\bul_Bd\bul_Bq$, i.e.
$\epsilon=\mu\circ(p\bul_B\varepsilon_d\bul_Bq)$.

To obtain a formula for $\zeta$, we calculate \equref{cotensor} for $n=x\bul_Bq$, where $(x,\rho)\in\Rcom(\Omega,\dd)$, then we find the following diagram.
\begin{equation}\eqlabel{zeta}
\xymatrix{
(x\bul_Bq)\bul^c(p\bul_Bd)\ar[r] &
(x\bul_Bq)\bul_A(p\bul_Bd)\ar@<-.5ex>[r] \ar@<.5ex>[r] &
(x\bul_Bq)\bul_Ac\bul_A(p\bul_Bd)\\
& x \ar[ul]^{\zeta_x} \ar[u]_{(x\ubul_B\eta\ubul_Bd)\circ\rho}
}
\end{equation}
We obtain $\zeta_x$ by the universal property of the equalizer.

\subsection{The canonical morphisms}

Consider $0$--cells $\Omega$,
$\Omega'$ and $B$, two $1$--cells $p\in\Hom_1(\Omega,B)$ and
$q\in\Hom_1(B,\Omega')$, and let
$\dd=(B,d,\Delta_d,\varepsilon_d)$ be a comonad in $\Bb$. Then
$(p\bul_Bd,p\ubul_B\Delta_d)\in\Rcom(\Omega,d)$ and
$(d\bul_Bq,\Delta_d\ubul_Bq)\in\Rcom(d,\Omega')$. With this
notation, we can state the following lemma.

\begin{lemma}\lelabel{can1}
The following isomorphism holds, in particular, the following equalizer (cotensor product) always exists in $\Hom(\Omega,\Omega')$,
$$(p\bul_Bd)\bul^d(d\bul_Bq)\cong p\bul_Bd\bul_Bq.$$
\end{lemma}

\begin{proof}
We claim that the equalizer of $p\ubul_B\Delta_d\ubul_Bd\ubul_Bq$ and $p\ubul_Bd\ubul_B\Delta_d\ubul_Bq$ is given by
$(p\bul_Bd\bul_Bq,p\ubul_B\Delta_d\ubul_Bq)$.
Consider any $2$--cell $\kappa\in{^\Omega\Hom_2^{\Omega'}}(x,p\ubul_Bd\ubul_Bd\ubul_Bq)$, such that $(p\ubul_B\Delta_d\ubul_Bd\ubul_Bq)\circ\kappa=(p\ubul_Bd\ubul_B\Delta_d\ubul_Bq)\circ\kappa$.
Then we find that
\begin{eqnarray*}
(p\ubul_B\Delta_d\ubul_Bq)\circ(p\ubul_Bd\ubul_B\varepsilon_d\ubul_Bq)\circ\kappa
&=& (p\ubul_Bd\ubul_Bd\ubul_B\varepsilon_d\ubul_Bq)\circ
(p\ubul_B\Delta_d\ubul_Bd\ubul_Bq)\circ\kappa\\
&=& (p\ubul_Bd\ubul_Bd\ubul_B\varepsilon_d\ubul_Bq)\circ
(p\ubul_Bd\ubul_B\Delta_d\ubul_Bq)\circ\kappa=\kappa,\\
\end{eqnarray*}
where we applied \equref{interchange} in the first equality and the condition on $\kappa$ in the second one.
Consequently, the following diagram commutes
\[
\xymatrix{
p\bul_Bd\bul_Bq \ar[rr]^-{p\ubul_B\Delta_d\ubul_Bq} && p\bul_Bd\bul_Bd\bul_Bq
\ar@<.5ex>[rr]^-{p\ubul_B\Delta_d\ubul_Bd\ubul_Bq} \ar@<-.5ex>[rr]_-{p\ubul_Bd\ubul_B\Delta_d\ubul_Bq} && p\bul_Bd\bul_Bd\bul_Bd\bul_Bq, \\
x \ar[u]^\lambda \ar[urr]_-{\kappa}
}
\]
where $\lambda=(p\ubul_Bd\ubul_B\varepsilon_d\ubul_Bq)\circ\kappa$.
Moreover, for any other $\lambda'\in{^\Omega\Hom_2^{\Omega'}}(x,p\ubul_Bd\ubul_Bq)$ such that $\kappa=(p\ubul_B\Delta_d\ubul_Bq)\circ\lambda'$, we find that
\begin{eqnarray*}
\lambda&=&(p\ubul_Bd\ubul_B\varepsilon_d\ubul_Bq)\circ\kappa\\
&=&(p\ubul_Bd\ubul_B\varepsilon_d\ubul_Bq)\circ(p\ubul_B\Delta_d\ubul_Bq)\circ\lambda'
=\lambda',
\end{eqnarray*}
i.e. $\lambda$ is the unique $2$--cell satisfying this property.
\end{proof}

Consider now a comonad-morphism with adjunction $(p,q):\dd\to\cc$.
Denote the associated adjoint pair by $(A,B,p,q,\eta,\mu)$.
We know that $p\bul_Bd\in\Bicom(\cc,\dd)$ and $d\bul_Bq\in\Bicom(\dd,\cc)$ (see \leref{bicomodule}).
With this notation the following lemma holds.

\begin{lemma}\lelabel{can2}
The map $\gamma=(d\ubul_B\eta\ubul_Bd)\circ\Delta_d:d\to d\bul_Bq\bul_Ap\bul_Bd$ satisfies the following equality
$$(d\ubul_Bq\ubul_A\lambda^{p\bul_Bd})\circ\gamma=(\rho^{d\bul_Bq}\ubul_Ap\ubul_Bd)\circ\gamma.$$
\end{lemma}

\begin{proof}
Let us denote by $(q,\alpha):\dd\to\cc$ the right comonad-morphism associated to the comonad-morphism with adjunction $(p,q)$. Then we know that
\begin{eqnarray*}
\rho^{d\bul_Bq}&=&(d\ubul_B\alpha)\circ(\Delta_d\ubul_Bq);\\
\lambda^{p\bul_Bd}
&=&(\mu\ubul_Ac\ubul_Ap\ubul_Bd)\circ(p\ubul_B\alpha\ubul_Ap\ubul_Bd)
\circ(p\ubul_Bd\ubul_B\eta\ubul_Bd)\circ(p\ubul_B\Delta_d).
\end{eqnarray*}
Then we can compute on one hand
\begin{eqnarray*}
(\rho^{d\bul_Bq}\ubul_Ap\ubul_Bd)\circ\gamma
&=&(d\ubul_B\alpha\ubul_Ap\ubul_Bd)\circ(\Delta_d\ubul_Bq\ubul_Ap\ubul_Bd)
\circ(d\ubul_B\eta\ubul_Bd)\circ\Delta_d\\
&=&(d\ubul_B\alpha\ubul_Ap\ubul_Bd)\circ(d\ubul_Bd\ubul_B\eta\ubul_Bd)\circ(\Delta_d\ubul_Bd)\circ\Delta_d,
\end{eqnarray*}
where we used \equref{interchange2} in the second equality. On the other hand, we find
\begin{eqnarray*}
(d\ubul_Bq\ubul_A\lambda^{p\bul_Bd})\circ\gamma
&=&(d\ubul_Bq\ubul_A\mu\ubul_Ac\ubul_Ap\ubul_Bd)
\circ(d\ubul_Bq\ubul_Ap\ubul_B\alpha\ubul_Ap\ubul_Bd)
\\
&&\hspace{.5cm}\circ(d\ubul_Bq\ubul_Ap\ubul_Bd\ubul_B\eta\ubul_Bd)
\circ(d\ubul_Bq\ubul_Ap\ubul_B\Delta_d)
\circ(d\ubul_B\eta\ubul_Bd)\circ\Delta_d\\
&=&(d\ubul_Bq\ubul_A\mu\ubul_Ac\ubul_Ap\ubul_Bd)
\circ(d\ubul_Bq\ubul_Ap\ubul_B\alpha\ubul_Ap\ubul_Bd)
\\
&&\hspace{.5cm}\circ(d\ubul_Bq\ubul_Ap\ubul_Bd\ubul_B\eta\ubul_Bd)
\circ(d\ubul_B\eta\ubul_Bd\ubul_Bd)
\circ(d\ubul_B\Delta_d)\circ\Delta_d\\
&=&(d\ubul_Bq\ubul_A\mu\ubul_Ac\ubul_Ap\ubul_Bd)
\circ(d\ubul_B\eta\ubul_Bq\ubul_Ac\ubul_Ap\ubul_Bd)
\circ(d\ubul_B\alpha\ubul_Ap\ubul_Bd)
\\
&&\hspace{.5cm}\circ(d\ubul_Bd\ubul_B\eta\ubul_Bd)
\circ(\Delta_d\ubul_Bd)\circ\Delta_d\\
&=&(d\ubul_B\alpha\ubul_Ap\ubul_Bd)
\circ(d\ubul_Bd\ubul_B\eta\ubul_Bd)
\circ(\Delta_d\ubul_Bd)\circ\Delta_d.\\
\end{eqnarray*}
Here we used \equref{interchange2} in the second and third equality, the coassociativity of $\dd$ as well in the third equality and \equref{ccc} in the last equality.
\end{proof}

\begin{Definition}
Consider a comonad-morphism with adjunction $(p,q):\dd\to\cc$.
Suppose that the cotensor product $(d\bul_Bq)\bul^c(p\bul_Bd)$ exists in $\Hom(B,B)$. The following canonical $2$--cell is well-defined by
\prref{locdual} and \leref{can1}
$$\ul{\can}:(p\bul_Bd)\bul^d(d\bul_Bq)\cong p\bul_Bd\bul_Bq \to c.$$
\leref{can2} together with the universal property of the equalizer, induces a well-defined $2$--cell
$$\ol{\can}:d \to (d\bul_Bq)\bul^c(p\bul_Bd).$$
If both $\ul{\can}$ and $\ol{\can}$ are isomorphisms, then we say
that $(p,q)$ is a $\dd$-$\cc$ \emph{Galois comonad-morphism}.
\end{Definition}

\subsection{Structure theorems}\selabel{structure}

We will now give necessary and sufficient conditions for the
pullback and pushout functor to be full and faithful, and then to
obtain an equivalence between the categories $\Rcom(\Omega,\cc)$
and $\Rcom(\Omega,\dd)$.

Consider the following diagram in $\Rcom(\Omega,\cc)$.
\begin{equation}\eqlabel{Psinu}
\xymatrix{
(y\bul^c(p\bul_Bd))\bul_Bq  \ar[d]_{\nu_y} \ar[r]^-{{\rm eq}_y\ubul_Bq}
& y\bul_A(p\bul_Bd)\bul_Bq \ar@<.5ex>[rr]^-{\rho^y\ubul_Ap\ubul_Bd\ubul_Bq} \ar@<-.5ex>[rr]_-{y\ubul_A\lambda^{p\bul_Bd}\ubul_Bq} \ar[d]^{y\ubul_A\ul{\can}} \ar[dl]^{y\ubul_A\epsilon}
&& y\bul_Ac\bul_Ap\bul_Bd\bul_Bq \ar[d]^{y\ubul_Ac\ubul_A\ul{\can}} \\
y \ar[r]_-{\rho^y}
& y\bul_Ac \ar@<.5ex>[rr]^-{\rho^y\ubul_Ac} \ar@<-.5ex>[rr]_-{y\ubul_A\Delta_c}
&& y\bul_Ac\bul_Ac
}
\end{equation}
As at the end of \seref{popb}, we have denoted $\epsilon=\mu\circ(p\ubul_B\varepsilon_d\ubul_Bq)$.

\begin{lemma}\lelabel{commute}
Consider the diagram \equref{Psinu} and any $\psi:x\to y\bul_Ap\bul_Bd\bul_Bq$ such that $(\rho^y\ubul_Ap\ubul_Bd\ubul_Bq)\circ \psi =(y\ubul_A\lambda^{p\bul_Bd}\ubul_Bq)\circ \psi$, i.e. $\psi$ equalizes the two upper horizontal arrows. Then
$$(y\ubul_A\ul{\can})\circ \psi=\rho^y\circ(y\ubul_A\epsilon)\circ \psi.$$
Consequently the left square of diagram \equref{Psinu} commutes and moreover the right square commutes in a serial way.
\end{lemma}

\begin{proof}
We compute
\begin{equation}\eqlabel{eqbar}
\begin{array}{rcl}
\rho^y\circ(y\ubul_A\epsilon)\circ \psi
&=&(y\ubul_Ac\ubul_A\epsilon)\circ(\rho^y\ubul_Ap\ubul_Bd\ubul_Bq)\circ \psi\\
&=&(y\ubul_Ac\ubul_A\epsilon)\circ(y\ubul_A\lambda^{p\bul_Bd}\ubul_Bq)\circ \psi\\
&=&(y\ubul_A\ul{\can})\circ \psi.
\end{array}
\end{equation}
Here we used \equref{interchange2} in the first equality, the equalizing condition of $\psi$ in the second equality and one of the alternative formulas for $\ul{\can}$ in the last equality.
We can conclude that the left square of \equref{Psinu} commutes by taking $(x,\psi)=((y\bul^c(p\bul_Bd))\bul_Bq,{\rm eq}_y\bul_Bq)$. The last statement is easy.
\end{proof}

The following theorem generalizes \cite[Theorem 5.1]{BrzElkGT} and
\cite[Theorem 2.7]{GT:comonad}. For its proof, we need
\leref{cotensor}, that states that ``all comodules are
equalizers'', which could have been gathered from Beck's Theorem.
However, we find it adequate to include a proof here.

\begin{lemma}\lelabel{cotensor}
Let $\cc=(A,c,\Delta_c,\varepsilon_c)$ be a comonad in $\Bb$ and
consider $(m,\rho)\in\Rcom(\Omega,\cc)$. Then the equalizer of
$(m\ubul_A\Delta_c)$ and $(\rho\ubul_Ac)$ always exists in
$\Rcom(\Omega,\cc)$ and is isomorphic to $(m,\rho)$, i.e.
$m\bul^cc\cong m$.
\end{lemma}
\begin{proof}
For any pair $(x,\xi)$ such that $\xi$ equalizes $m\ubul_A \Delta$
and $\rho\ubul_A c$, we can define
$\zeta=(m\ubul_A\varepsilon)\circ\xi$. We have that
$\rho\circ\zeta=\xi$, indeed:
\begin{gather}
\rho\circ\zeta=\rho\circ(m\ubul_A\varepsilon)\circ\xi =
(m\ubul_Ac\ubul_A\varepsilon)\circ(\rho\ubul_Ac)\circ\xi =
(m\ubul_Ac\ubul_A\varepsilon)\circ(m\ubul_A\Delta)\circ\xi = \xi
\end{gather}
Moreover if there exists a $\zeta'$ such that
$\rho\circ\zeta'=\xi$, then we find
$(m\ubul_A\varepsilon)\circ\xi=(m\ubul_A\varepsilon)\circ\rho\circ\zeta'=\zeta'$
\end{proof}

\begin{theorem}\thlabel{weakstructure}
Let $\dd$, $\cc$ be comonads in $\Bb$ and consider a
comonad-morphism with adjunction $(p,q):\dd\to\cc$ (see
\prref{locdual}). Suppose $\Rcom(\Omega,\dd)$ satisfies the
equalizer condition for $p$. Then the pullback functor $\pP$
associated to $p$ is fully faithful if and only if $\ul{\can}$ is
an isomorphism and the pushout functor $\qQ$ preserves the
equalizers of the form \equref{cotensor}.
\end{theorem}

\begin{proof}
Suppose first that $\ul{\can}$ is an isomorphism and $\qQ$
preserves the equalizers of the form \equref{cotensor}. This last
condition means that $(y\bul^c(p\bul_Bd))\bul_Bq\cong
y\bul^c(p\bul_Bd\bul_Bq)$. Moreover, since $\ul{\can}$ is an
isomorphism, we can compute further $y\bul^c(p\bul_Bd\bul_Bq)\cong
y\bul^cc$, and by \leref{cotensor} we know that the equalizer
$y\bul^cc=y$, so we conclude $(y\bul^c(p\bul_Bd))\bul_Bq\cong y$
and the pullback functor is fully faithful.

To prove the converse, take first $y=c$, then we see that the equalizer $(c\bul^c(p\bul_Bd),{\rm eq}_c)\cong (p\bul_Bd,\lambda^{p\bul_Bd})$ by \leref{cotensor}, and consequently $(c\bul^c(p\bul_Bd))\bul_Bq\cong p\bul_Bd\bul_Bq$. As a consequence, we find that
\begin{eqnarray*}
\nu_c&=&(c\bul_A\epsilon)\circ({\rm eq}_c\bul_Bq)\\
&=&(c\bul_A\epsilon)\circ(\lambda^{p\bul_Bd}\bul_Bq)\\
&=&(c\bul_A\epsilon)\circ(\beta\bul_Bd\bul_Bq)\circ(p\bul_B\Delta_d\bul_Bq)\\
&=&(c\bul_A\epsilon)\circ(\ul{\can}\bul_Ap\bul_Bd\bul_Bq)\circ(p\bul_Bd\bul_B\eta\bul_Bd\bul_Bq)\circ(p\bul_B\Delta_d\bul_Bq)\\
&=&\ul{\can}\circ(p\bul_Bd\bul_Bq\bul_B\epsilon)\circ\Delta_{p\bul_Bd\bul_Bq}
=\ul{\can}.
\end{eqnarray*}
If the pullback functor is full and faithful, then $\nu$ is a
natural isomorphism, so in particular, we find that
$\nu_c=\ul{\can}$ is an isomorphism.

Take any $(y,\rho^y)$ in $\Rcom(\Omega,\cc)$. We are done if we
show that $((y\bul^c(p\bul_Bd))\bul_Bq,{\rm eq}_y\ubul_Bq)$ is the
equalizer of $\rho^y\ubul_Ap\ubul_Bd\ubul_Bq$ and
$y\ubul_A\lambda^{p\bul_Bd}\ubul_Bq$. To this end, consider the
diagram \equref{Psinu}, where all squares are (serial) commutative
by \leref{commute}. Since we know that all vertical lines are
isomorphisms and the lower horizontal line is an equalizer, we
find that the upper horizontal line must be an equalizer as well.
\end{proof}

In order to characterize when the pullback functor is fully faithful, we need the following lemma.

\begin{lemma}\lelabel{equalizerxq}
Let $\dd=(B,d,\Delta_d,\varepsilon_d)$ be a comonad in $\Bb$ and consider $(x,\rho)\in \Rcom(\Omega,\dd)$ and $q\in \Hom_1(B,A)$. Then
$(x\bul_Bq,\rho\bul_Bq)$ is the equalizer of $\rho\bul_Bd\bul_Bq$ and $x\bul_B\Delta_d\bul_Bq$ in $\Hom(\Omega,A)$, i.e. we have the equalizer
\begin{equation}\eqlabel{equalizerxq}
\xymatrix{
x\bul_Bq\cong x\bul^d(d\bul_Bq) \ar[rr] && x\bul_Bd\bul_Bq \ar@<.5ex>[rr]^-{\rho\bul_Bd\bul_Bq}
\ar@<-.5ex>[rr]_-{x\bul_B\Delta_d\bul_Bq} && x\bul_Bd\bul_Bd\bul_Bq.
}
\end{equation}
Moreover, if $\cc=(A,c,\Delta_c,\varepsilon_c)$ is another comonad in $\Bb$ and
$d\bul_Bq$ is a $\dd$-$\cc$--bicomodule with left $\dd$--coaction $\Delta_d\bul_Bq$,
then the above equalizer exists even in $\Rcom(\Omega,\cc)$.
\end{lemma}

\begin{proof}
The first part of the lemma can be proven in the same way as \leref{can1}.

Suppose that
$d\bul_Bq$ is a $\dd$-$\cc$--bicomodule with left $\dd$--coaction $\Delta_d\bul_Bq$, and right $\cc$--coaction $\rho'$.
Then one can verify that $(x\bul_Bq,\varrho)\in\Rcom(\Omega,\cc)$, where $\varrho=(x\ubul_B\varepsilon_d\ubul_Bq\ubul_Ac)\circ(x\ubul_B\rho')\circ(\rho\ubul_Bq)$.
Furthermore, one can use the same technique as in the proof of \leref{equalizer} to check that the equalizer \equref{equalizerxq} exists in $\Rcom(\Omega,\cc)$.
\end{proof}

\begin{theorem}\thlabel{Qfullyfaithful}
Let $\dd$, $\cc$ be comonads in $\Bb$ and consider a comonad-morphism
with adjunction $(p,q):\dd\to\cc$ (see \prref{locdual}).
 Suppose $\Rcom(\Omega,\dd)$ satisfies the equalizer condition
for $p$ and that the equalizer $(d\bul_Bq)\bul^c(p\bul_Bd)$ is
left $d\bul_Bd$--pure. Then the following statements are
equivalent:
\begin{enumerate}[(i)]
\item
The pullback functor $\qQ$ associated to $q$ is fully faithful;
\item $\ol{\can}$ is an isomorphism and
\begin{equation}\eqlabel{cotensorass}
(x\bul_Bq)\bul^c(p\bul_Bd)\cong x\bul^d((d\bul_Bq)\bul^c(p\bul_Bd))
\end{equation}
for all $x\in\Rcom(\Omega,\dd)$;
\item $\ol{\can}$ is an isomorphism and
the equalizer $(d\bul_Bq)\bul^c(p\bul_Bd)$ is left $y$--pure for all $y\in\Hom_1(\Omega,B)$.
\end{enumerate}
\end{theorem}

\begin{proof}
$\ul{(i)\Leftrightarrow(ii)}$.
Since the equalizer $(d\bul_Bq)\bul^c(p\bul_Bd)$ is left $d\bul_Bd$--pure, $(d\bul_Bq)\bul^c(p\bul_Bd)$ is a left $\dd$--comodule and the equalizer $x\bul^d((d\bul_Bq)\bul^c(p\bul_Bd))$ is well-defined for all $x\in\Rcom(\Omega,\dd)$.

It follows directly from the definition of $\ol{\can}$ and the formula for $\zeta$ in \equref{zeta} that $\zeta_d=\ol{\can}$.
If $\ol{\can}$ is an isomorphism, than we find
$$x\bul^d((d\bul_Bq)\bul^c(p\bul_Bd))\cong x\bul^d d\cong x.$$
Combining this isomorphism with \equref{cotensorass}, we find that
$$\pP\qQ(x)=(x\bul_Bq)\bul^c(p\bul_Bd)\cong x,$$
and from naturalness it follows that this combined isomorphism is
exactly the counit of the adjunction evaluated in $x$.
The converse follows in the same way.\\
$\ul{(ii)\Rightarrow(iii)}$. Take any $y\in\Hom_1(\Omega,B)$ and put $x=y\bul_Bd$ in \equref{cotensorass}, then we find
$$(y\bul_Bd\bul_Bq)\bul^c(p\bul_Bd)\cong
(y\bul_Bd)\bul^d((d\bul_Bq)\bul^c(p\bul_Bd))\cong y\bul_B((d\bul_Bq)\bul^c(p\bul_Bd)),
$$
where we used \leref{equalizerxq} for the last isomorphism. This means exactly that $(d\bul_Bq)\bul^c(p\bul_Bd)$ is left $y$--pure.\\
$\ul{(iii)\Rightarrow(ii)}$. If $(d\bul_Bq)\bul^c(p\bul_Bd)$ is
left $y$--pure for all $y\in\Hom_1(\Omega,B)$, then
$$y\bul_B((d\bul_Bq)\bul^c(p\bul_Bd))\cong(y\bul_Bd\bul_Bq)\bul^c(p\bul_Bd).$$
Now take $x\in\Rcom(\Omega,\dd)$ and construct the following
commutative diagram.
\[
\hspace{-.11cm}
\xymatrix{
(x\bul_Bq)\bul^c(p\bul_Bd) \ar[r] & (x\bul_Bd\bul_Bq)\bul^c(p\bul_Bd) \ar@<.5ex>[r] \ar@<-.5ex>[r] & (x\bul_Bd\bul_Bd\bul_Bq)\bul^c(p\bul_Bd) \\
x\bul^d((d\bul_Bq)\bul^c(p\bul_Bd)) \ar[r] \ar[u]
& x\bul_B((d\bul_Bq)\bul^c(p\bul_Bd)) \ar@<.5ex>[r] \ar@<-.5ex>[r] \ar[u]_\cong
& x\bul_Bd\bul_B((d\bul_Bq)\bul^c(p\bul_Bd)) \ar[u]_\cong
}
\]
The upper row is obtained by applying the functor $\pP$ on the equalizer \equref{equalizerxq}. Since $\pP$ has a left adjoint, it preserves the equalizers and the upper row is an equalizer. The lower row is an equalizer by definition. The vertical isomorphisms follow from our previous observation. By the universal property of the equalizer we obtain that the left vertical arrow must be an isomorphism as well.
\end{proof}

\begin{theorem}\thlabel{strongstructure}
Let $\dd$, $\cc$ be comonads in $\Bb$ and consider a
comonad-morphism with adjunction $(p,q):\dd\to\cc$ (see
\prref{locdual}). Suppose $\Rcom(\Omega,\dd)$ satisfies the
equalizer condition for $p$. Then the following statements are
equivalent.
\begin{enumerate}[(i)]
\item The functors $(\pP,\qQ)$ establish an equivalence of categories between $\Rcom(\Omega,\cc)$ and $\Rcom(\Omega,\dd)$;
\item $\ul{\can}$ is an isomorphism, the pushout functor $\qQ$ reflects isomorphisms and preserves the equalizers of the form \equref{cotensor}.
\item[]\hspace{-1.1cm}
\begin{minipage}{\textwidth}
\vspace{.1cm}
If $(d\bul_Bq)\bul^c(p\bul_Bd)$ is left $d\bul_Bd$--pure, then the previous statements are furthermore equivalent to
\vspace{.2cm}
\end{minipage}
\item $(p,q)$ is a Galois comonad-morphism, $\qQ$ preserves the equalizers of the form \equref{cotensor} and \equref{cotensorass} is satisfied for all $x\in\Rcom(\Omega,\dd)$;
\item $(p,q)$ is a Galois comonad-morphism, $\qQ$ preserves the equalizers of the form \equref{cotensor} and $(d\bul_Bq)\bul^c(p\bul_Bd)$ is left $y$--pure for all $y\in\Hom_1(\Omega,B)$.
\end{enumerate}
\end{theorem}

\begin{proof}
We only prove the equivalence between (i) and (ii), the equivalence
with the other statements follows directly from
\thref{weakstructure} and \thref{Qfullyfaithful}.

First suppose that $(\pP,\qQ)$ establishes an equivalence of
categories. Then obviously $\qQ$ reflects isomorphisms, and the
other statements follow from  \thref{weakstructure}.

Conversely, suppose $\ul{\can}:p\bul_Bd\bul_Bq \rightarrow c$ to
be an isomorphism of comonads, then by \leref{cotensor}, we find
the following equalizer in $\Rcom(\Omega,\cc)$ for any
$(m,\rho^m)$ in $\Rcom(\Omega,\dd)$,
$$
\xymatrix{ m\bul_Bq \ar[r] & m\bul_Bq\bul_Ap\bul_Bd\bul_Bq
\ar@<.5ex>[r] \ar@<-.5ex>[r] &
m\bul_Bq\bul_Ap\bul_Bd\bul_Bq\bul_Ap\bul_Bd\bul_Bq. }$$
Furthermore, if $\qQ$ preserves equalizers, then we can apply
$-\bul_Bq$ on \equref{cotensor} in the situation $n=m\bul_Bq$ and
we obtain a second equalizer in $\Rcom(\Omega,\cc)$. These two
equalizers can be related in the following diagram.
\[
\xymatrix{
 m\bul_Bq \ar[rr]^-{\zeta_m\bul_Bq} \ar[d]
&& ((m\bul_Bq)\bul^c(p\bul_Bd))\bul_Bq \ar[d]^{{\rm eq}\ubul_Bq}\\
 (m\bul_Bq)\bul_A(p\bul_Bd\bul_Bq) \ar@<.5ex>[d] \ar@<-.5ex>[d] \ar@{=}[rr]
&&(m\bul_Bq)\bul_A(p\bul_Bd)\bul_Bq \ar@<.5ex>[d] \ar@<-.5ex>[d] \\
 (m\bul_Bq)\bul_A(p\bul_Bd\bul_Bq)\bul_A(p\bul_Bd\bul_Bq) \ar[rr]^-{m\ubul_Bq\ubul_A\ul{\can}\ubul_Ap\ubul_Bd\ubul_Bq}
&&(m\bul_Bq)\bul_Ac\bul_A(p\bul_Bd)\bul_Bq
}
\]
Since $\ul{\can}$ is an isomorphism, we find by the properties of the equalizers that $\zeta_m\bul_Bq$ is an isomorphism as well, and since $\qQ$ reflects isomorphisms, $\zeta_m$ must be an isomorphism. From \thref{weakstructure} we know that $\nu$ is a natural isomorphism as well, so we find that $(\pP,\qQ)$ is an equivalence of categories.
\end{proof}

\subsection{Coseparable comonads}

A coalgebra in a monoidal category is called coseparable if its
comultiplication splits in the category of bicomodules. In
analogy, we will say that a comonad
$\cc=(A,c,\Delta_c,\varepsilon_c)$ is coseparable if and only if
there exists a $2$--cell $\gamma\in {^A\Hom^A_2}(c\bul_Ac,c)$ which
is a morphism of $\cc$-$\cc$--bicomodules and such that
$\gamma\circ\Delta_c=c$.

The following proposition has a straightforward proof.

\begin{proposition}\prlabel{cosepcomonad1}
Consider a quadruple $\cc=(A,c,\Delta,\varepsilon)$ in $\Bb$
consisting of a $0$--cell $A$, a $1$--cell $c:A\to A$ and $2$--cells
$\Delta : c\to c\bul c$ and $\varepsilon: A\to c$. The following
statements are equivalent.
\begin{enumerate}[(i)]
\item $\cc$ is a coseparable comonad in $\Bb$; \item
$(c,\Delta,\varepsilon)$ is a coseparable coalgebra in the
monoidal category $\Hom(A,A)$; \item $(- \bul_A c, - \ubul_A
\Delta,\ubul_A \varepsilon)$ is a coseparable comonad on
$\Hom(\Omega,A)$ for all $\Omega\in\Bb$; \item $(- \bul_A c, -
\ubul_A \Delta,\ubul_A \varepsilon)$ is a coseparable comonad on
$\Hom(A,A)$.
\end{enumerate}
\end{proposition}
The notion of a separable functor was introduced in \cite{NVV} and
can be used to characterize separable algebras and coseparable
coalgebras (see e.g. \cite{CMZ}), and, more generally, of
separable morphisms of corings \cite{Gomez:2002}. Separable
functors of the second kind were introduced in \cite{CM}. A useful
characterization of separable functors is given by Rafael's
theorem \cite{rafael}. Let us state a version of Rafael's theorem
for separable functors of the second kind (the corresponding
theorem for separable functors of the first kind can be obtained
by taking $\Cc=\Aa$ and $H=\id_\Aa$).

\begin{theorem}[{\cite[Theorem 2.7]{CM}}] Let $F:\Aa\to\Bb$ be a covariant functor with a right adjoint $G:\Bb\to\Aa$ and consider an additional functor $H:\Aa\to\Cc$. Denote by $\eta:\id_\Aa\to GF$ the unit of the adjunction $(F,G)$, then $F$ is $H$--separable if and only if there exists a natural transformation $\mu:HGF\to H$ such that $\mu\circ H\eta=H$.
\end{theorem}

Combining \prref{cosepcomonad1} with \cite[Theorem 1.6]{EKV}, we
obtain more characterizations of coseparable comonads in
bicategories.
\begin{theorem}
Consider a comonad $\cc=(A,c,\Delta,\varepsilon)$ in $\Bb$. The
following statements are equivalent.
\begin{enumerate}[(i)]
\item $\cc$ is coseparable; \item the forgetful functor
$U_{A,\cc}:\Rcom(A,\cc)\to\Hom(A,A)$ is separable; \item the
forgetful functor
$U_{\Omega,\cc}:\Rcom(\Omega,\cc)\to\Hom(\Omega,A)$ is separable
for all $0$--cells $\Omega$ in $\Bb$; \item the forgetful functor
${_{\cc,A}U}:\Lcom(\cc,A)\to\Hom(A,A)$ is separable; \item the
forgetful functor
${_{\cc,\Omega}U}:\Lcom(\cc,\Omega)\to\Hom(A,\Omega)$ is separable
for all $0$--cells $\Omega$ in $\Bb$.
\end{enumerate}
\end{theorem}
\begin{remark}\relabel{charcosep}
Recall that the unit of the adjunction $U_{\Omega,\cc} \dashv
-\bul_Ac$ is given by
$\theta_m=\rho^m\in{^\Omega\Hom^A_2}(m,m\bul_Ac)$, the right
$\cc$--coaction for all $(m,\rho^m)\in \Rcom(\Omega,\cc)$.
Therefore, a comonad $\cc$ is coseparable if and only if there
exists, for all $0$--cells $\Omega$ in $\Bb$ and all $(m,\rho^m)\in
\Rcom(\Omega,\cc)$, a $2$--cell
$\bar{\rho}_m\in{^\Omega\Hom^A_2}(m\bul_Ac,m)$ such that
$\bar{\rho}_m$ is natural in $m$ and $\bar{\rho}_m\circ\rho^m=m$.
This is equivalent to the existence, for all $0$--cell $\Omega$ in
$\Bb$ and all $(n,\lambda^n)\in \Lcom(\cc,\Omega)$, of a $2$--cell
$\bar{\lambda}_n\in{^A\Hom^\Omega_2}(c\bul_An,n)$ such that
$\bar{\lambda}_n$ is natural in $n$ and
$\bar{\lambda}_n\circ\lambda^n=n$.~\qed
\end{remark}

Our first aim is to prove that Galois comonad-morphisms having as
codomain a coseparable comonad give rise to equivalences of
categories. The following lemma generalizes
\cite[2.13]{Wis:galcom} and \cite[Lemma 9.1]{V:PhD}.

\begin{lemma}\lelabel{cotensorinvcoproduct}
Let $\cc=(A,c,\Delta_c,\varepsilon_c)$ be a comonad in $\Bb$ and
$0$--cells $B$, $\Omega$ and $\Omega'$ in $\Bb$. Take $(m,\rho^m)
\in \Rcom(\Omega,\cc)$ and $(n,\lambda^n) \in \Lcom(\cc,B)$.
Suppose that the cotensor product $m\bul^cn$ exists in
$\Hom(\Omega,B)$, then for all $p\in \Hom(B,\Omega')$ the functor
$- \bul_Ap:\Hom(\Omega,A)\to\Hom(\Omega,B)$ preseves this
equalizer, i.e. natural map
$$f:\ (m\bul^cn)\bul_B p\to m\bul^c (n\bul_Bp)$$
is an isomorphism in each of the following situations:
\begin{enumerate}[(i)]
\item $\rho^m:m\to m\bul_Ac$ has a left inverse $\bar{\rho}_m$ in $\Rcom(\Omega,\cc)$;
\item $\lambda^n:n\to c\bul_An$ has a left inverse $\bar{\lambda}_n$ in $\Lcom(\cc,B)$.
\end{enumerate}
\end{lemma}

\begin{proof}
Let us show that the equalizer defining the cotensor product
$$\xymatrix{
m\bul^c n \ar[r]^i & m\bul_A n \ar@<.5ex>[rr]^-{\rho^m\bul_An}
\ar@<-.5ex>[rr]_-{m\bul_A\lambda^n} && m\bul_Ac\bul_An },$$ is a
contractible equalizer. We define $2$--cells $\alpha : m\bul_A n\to
m\bul^c n$ and $\beta,\beta':m\bul_Ac\bul_An\to m\bul_An$ as
follows. Under condition (i) we put $\alpha=
(\bar{\rho}_m\ubul_An)\circ(m\ubul_A\lambda^n)$ and
$\beta=\bar{\rho}_m\ubul_An$, if (ii) is satisfied, then we define
$\alpha= (m\ubul_A\bar{\lambda}_n)\circ(\rho^m\ubul_An)$ and
$\beta'=m\ubul_A\bar{\lambda}_n$.  Moreover, one can easily verify
that $\alpha\circ i= m\bul^c n$, $\beta\circ
(\rho^m\ubul_An)=m\bul_An=\beta'\circ (m\ubul_A\lambda^n)$ and
$\alpha\circ
i=\beta\circ(n\ubul_A\lambda^n)=\beta'\circ(\rho^m\ubul_An)$. By
\cite[Proposition 3.3.2]{BarrWells:ttt}, any functor preserves a
contractible equalizer. Therefore, we obtain the following
diagram, applying the functor $-\bul_Bp:\Hom(\Omega,B)\to
\Hom(\Omega,\Omega')$.
$$\xymatrix{
(m\bul^{c} n)\bul_B p\ar[r] \ar[d]^f &m\bul_A n\bul_B p \ar@<.5ex>[r] \ar[d]^\cong \ar@<-.5ex>[r] &
m\bul_A {c}\bul_A n\bul_B p \ar[d]^\cong \\
m\bul^{c} (n\bul_B p) \ar[r] & m\bul_A n\bul_B p \ar@<.5ex>[r] \ar@<-.5ex>[r] &
m\bul_A {c}\bul_A n\bul_B p}$$
Since we know that both horizontal rows are equalizers, $f$ is an isomorphism by the universal property of the equalizer.
\end{proof}

The following theorem generalizes \cite[Theorem 9.2]{V:PhD} and
\cite[5.7]{Wis:galcom}, it should be compared with \cite[Theorem
5.8]{CDV}. Our result is as well related to the split
(co)monadicity theorem (see e.g. \cite[Theorem
2.2]{JanTho:descIII}).

\begin{theorem}\thlabel{strcosep1}
Let $(p,q):\dd\to\cc$ be a comonad-morphism with adjunction and
suppose that $\Rcom(\Omega,\dd)$ satisfies the equalizer condition
for $p$. If $\lambda^{p\bul_Bd}: p\bul_Bd\to c\bul_Ap\bul_Bd$ has
a left inverse in $\Lcom(\cc,B)$ (in particular, if $\cc$ is a coseparable comonad), then the following statements hold.
\begin{enumerate}[(i)]
\item If $\ul{\can}$ is an isomorphism, then $\pP$ is fully
faithful;
\item if $\ol{\can}$ is an isomorphism, then $\qQ$ is
fully faithful;
\item if $(p,q)$ is a Galois comonad-morphism,
then $(\pP,\qQ)$ is an equivalence of categories.
\end{enumerate}
\end{theorem}

\begin{proof}
\ul{(i)}.
Since the coaction $\lambda^{p\bul_Bd} : p\bul_Bd\to c\bul_Ap\bul_Bd$ has a left inverse in $\Lcom(\cc,B)$, it follows by \leref{cotensorinvcoproduct} that for all $m\in\Rcom(\Omega,\cc)$, the following isomorphism holds
$$(m\bul^c(p\bul_Bd))\bul_Bq \cong m\bul^c((p\bul_Bd)\bul_Bq),$$
i.e. $\qQ$ preserves equalizers of the form \equref{cotensor}. The statement follows now by \thref{weakstructure}.\\
\ul{(ii)}.
Applying now left-right duality on \leref{cotensorinvcoproduct}, we find now that $(d\bul_Bq)\bul^c(p\bul_Bd)$ is left $y$--pure for all choices of $y\in\Hom_1(\Omega,B)$, since $\lambda^{p\bul_Bd}$ has a left inverse in $\Lcom(\cc,B)$. Therefore, $\qQ$ is fully faithful by \thref{Qfullyfaithful}.\\
\ul{(iii)}.
Follows directly from (i) and (ii).
\end{proof}

\begin{lemma}\lelabel{separable}
Let $F \dashv G$ and $H \dashv K$ be two pairs of adjoint functors
as in the following diagram
\[
\xymatrix{
\Aa \ar@<.5ex>[r]^F & \Bb \ar[d]^Z \ar@<.5ex>[l]^G \ar@<.5ex>[r]^H & \Cc \ar@<.5ex>[l]^K\\
 & \Dd
}
\]
Then the following statements hold.
\begin{enumerate}[(i)]
\item If $F$ is separable and $H$ is $G$--separable, then $HF$ is separable;
\item if $H$ is $Z$--separable, then $HF$ is $ZF$ separable.
\end{enumerate}
\end{lemma}

\begin{proof}
Denote the unit of the adjunction $F \dashv G$ by $\eta_a:a\to
GFa$ and the unit of $H \dashv K$ by $\zeta_b:b\to KHb$. Then the
unit of the composed adjoint pair $GK \dashv HF$ is given by
$$\xi_a=G\zeta_{Fa}\circ \eta_a : a\to GKHFa.$$
\ul{(i)}. If $F$ is separable, then there exists a natural transformation $\mu_a:GFa\to a$ which is a left inverse for $\eta_a$. If $H$ is $G$--separable, then there exists a natural transformation $\nu_b:GKHb\to Gb$ which is a left inverse for $G\zeta_b$. One can now easily see that $\mu_a\circ \nu_{Fa}$ is a left inverse for $\xi_a$.\\
\ul{(ii)}. Let $\varepsilon_b:FGb\to b$ be the counit of the adjunction $(F,G)$, and denote by $\nu_b:ZKHb\to Zb$ the inverse for $Z\zeta_b$, obtained by the $Z$--separability of $H$. Then we can consider the following diagram.
\[
\xymatrix{
ZFa \ar[r]^-{ZF\eta_a} \ar@{=}[rd] & ZFGFa \ar[r]^-{ZFG\zeta_{Fa}} \ar[d]^{Z\varepsilon_{Fa}} &
ZFGKHFa \ar[d]^{Z\varepsilon_{KHFa}}\\
& ZFa \ar[r]^-{Z\zeta_{Fa}} \ar@{=}[dr] & ZKHFa \ar[d]^{\nu_{Fa}}\\
&& ZFa
}
\]
The upper left triangle in this diagram commutes by the fact that
$F \dashv G$ is an adjunction, the lower right triangle commutes
by the fact that $H$ is $Z$--separable, the inner square commutes
trivially. Therefore we find that $ZF\xi_a=ZFG\zeta_{Fa}\circ
ZF\eta_a$ has a left inverse given by $\nu_{Fa}\circ
Z\varepsilon_{KHFa}$, i.e. $ZH$ is $ZF$--separable.
\end{proof}

\begin{remark}
It follows by \leref{separable} (ii) that any functor $F$ with a right adjoint is $F$--separable.~\qed
\end{remark}

The next Theorem is also related to the split (co)monadicity theorem (see e.g. \cite[Theorem 2.2]{JanTho:descIII}).

\begin{theorem}\thlabel{strcosep2}
Let $(p,q):\dd\to\cc$ be a comonad-morphism with adjunction, and
suppose that $\Rcom(\Omega,\dd)$ satisfies the equalizer condition
for $p$. Consider the forgetful
functor $U_{\Omega,\dd} : \Rcom(\Omega,\dd) \rightarrow
\Hom(\Omega,B)$ and its right adjoint
$G = - \bul_B d :\Hom(\Omega,B) \rightarrow
\Rcom(\Omega,\dd)$.

If $\ul{\can}$ is an isomorphism and
the functor $(- \bul_B q) U_{\Omega,\dd}$ is separable
(in particular, if the functor $- \bul_B q$ is $G$--separable and $\dd$ is
coseparable),
then $\pP$ is fully faithful.

If moreover $\ol{\can}$ is an isomorphism, then $(\pP,\qQ)$ is an equivalence of categories.
\end{theorem}

\begin{proof}
Obviously, $G (- \bul_A p)$ is right adjoint to $(-
\bul_B q)U_{\Omega,\dd}$. The unit is of this adjunction is given
by $\chi_m=(m\ubul_B\eta\ubul_Bd)\circ\rho^m:m\to
m\bul_Bq\bul_Ap\bul_Bd$ for all $(m,\rho^m)\in\Rcom(\Omega,\dd)$.
If the functor $(- \bul_B q) U_{\Omega,\dd}$ is separable, then for all $m\in \Hom(\Omega,\dd)$, the map $\chi_m:m\to m\bul_Bq\bul_Ap\bul_Bd$ has a left inverse. Taking $m=p\bul_Bd$ and taking into account the fact that $\ul{\can}$ is an isomorphism, we find that $\lambda^{p\bul_Bd}$ has a left inverse. Therefore, the statement follows by \thref{strcosep1}.\\
It Follows by \leref{separable} (i) that the functor $(- \bul_B q) U_{\Omega,\dd}$ is separable if
the functor $- \bul_B q$ is $G$--separable and $\dd$ is coseparable.
\end{proof}

The following theorem generalizes \thref{strcosep2} (taking $\Bb=\Rcom(\Omega,\dd)$, $\Yy=\Hom(\Omega,A)$, $\qQ=-\bul_Bq$ and $H$ and $H'$ the identity functors), however the proof of \thref{Joyal} is a lot more involved. \thref{Joyal} generalizes as well \cite[Theorem 2.3]{JanTho:descIII} and goes back to the Joyal-Tierney theorem for descent theory (see \seref{descent}).

\begin{theorem}\thlabel{Joyal}
Let $\cc=(A,c,\Delta_c,\varepsilon_c)$ and
$\dd=(B,d,\Delta_d,\varepsilon_d)$ be comonads in $\Bb$ and
$(p,q):\dd\to\cc$ be a comonad-morphism with adjunction such that
$\Rcom(\Omega,\dd)$ satisfies the equalizer condition for $p$.
Suppose furthermore that there exist functors $H$, $H'$ and $\qQ'$ that make the following diagram commutative
\[
\xymatrix{
\Rcom(\Omega,\dd) \ar[rr]^-\qQ \ar[dd]_H \ar[dr]_{U_{\Omega,\dd}} && \Rcom(\Omega,\cc) \ar[d]^{U_{\Omega,\cc}} \\
& \Hom(\Omega,B) \ar[r]^-{- \bul_B q} & \Hom(\Omega,A) \ar[d]^{H'}\\
\Bb \ar[rr]_{\qQ'} && \Yy
}
\]
(remark that the upper inner square of this diagram always
commutes, see \equref{diagramGal})
\begin{enumerate}
\item The functor $\pP$ is fully faithful if the following conditions hold,
\begin{enumerate}[(i)]
\item $\ul{\can}$ is an isomorphism
\item the composite functor $(- \bul_B q) U_{\Omega,\dd}$ is
$H$--separable; \item $H$ preserves equalizers of the form
\equref{cotensor}; \item $H'$ reflects isomorphisms.
\end{enumerate}
\item
The functor $\qQ$ is fully faithful and therefore $(\qQ,\pP)$ is an equivalence of categories if in addition to (1) the following conditions hold,
\begin{enumerate}[(i)]
\item $\ol{\can}$ is an isomorphism (hence $(p,q)$ is a Galois $\dd$-$\cc$ comonad-morphism);
\item $H$ reflects isomorphisms.
\end{enumerate}
\end{enumerate}
\end{theorem}

\begin{proof}
Denote, as in the proof of \thref{strcosep2}, the unit of the
adjoint pair
$$(- \bul_B q) U_{\Omega,\dd} \dashv (-\bul_Bd) (- \bul_A p)$$ by $$\chi_m=(m\ubul_B\eta\ubul_Bd)\circ\rho^m:m\to m\bul_Bq\bul_Ap\bul_Bd$$
for all $(m,\rho^m)\in\Rcom(\Omega,\dd)$. Then $(- \bul_B q)
U_{\Omega,\dd}$ is $H$--separable if and only if there exist
morphisms
$\xi_m:H(m\bul_Bq\bul_Ap\bul_Bd)\to H(m)$ for all $(m,\rho^m)\in\Rcom(\Omega,\dd)$ such that $\xi$ is natural in $m$ and $\xi_m\circ H(\chi_m)= H(m)$.\\
\ul{(1)}. By \thref{weakstructure}, we only have to prove that
$\qQ$ preserves equalizers of the form \equref{cotensor}. Since
$\ul{\can}$ is an isomorphism, we will identify $\cc$ with the
comonad $(p\bul_Bd\bul_Bq)$, the same applies for the
$\cc$--comodules. Consider the following commutative diagram of
equalizers for any $x\in\Rcom(\Omega,\cc)$.
\[
\xymatrix{
x\bul^c(p\bul_Bd) \ar[rr]^{\rm eq_x} \ar[d]_{\rm eq_x} && x\bul_A(p\bul_Bd) \ar[d]^{\rho^x\ubul_Ap\ubul_Bd}\\
x\bul_A(p\bul_Ad) \ar@<.5ex>[d]^{\chi_{x\bul_Ap\bul_Bd}} \ar@<-.5ex>[d]_{\rho^x\ubul_Ap\ubul_Bd} \ar[rr]^{\chi_{x\bul_Ap\bul_Bd}}
&& x\bul_Ac\bul_A(p\bul_Bd) \ar@<.5ex>[d]^{\chi_{x\bul_Ap\bul_Bd}\ubul_Bq\ubul_Ap\ubul_Bd} \ar@<-.5ex>[d]_{\rho^x\ubul_Ac\ubul_Ap\ubul_Bd} \\
x\bul_Ac\bul_A(p\bul_Bd) \ar[rr]^-{\chi_{x\bul_Ac\bul_Ap\bul_Bd}}
&& x\bul_Ac\bul_Ac\bul_A(p\bul_Bd)
}
\]
 The equalizers in both columns are of the form \equref{cotensor}, taking $n=x$
and $n=x\bul_Ac$ respectively. Consequently, if we apply the
functor $H$ to this diagram, then we obtain again a commutative
diagram, with equalizers in the columns. Moreover, we find that
the two lower horizontal arrows are split by respectively
$\xi_{x\bul_Ac\bul_Ap\bul_Bd}$ and $\xi_{x\bul_Ap\bul_Bd}$.
\begin{equation}
\eqlabel{equalizerdiagram}
\xymatrix{
H(x\bul^c(p\bul_Bd)) \ar@<.5ex>[rrr]^{H({\rm eq}_x)} \ar[d]_{H({\rm eq}_x)} &&& H(x\bul_A(p\bul_Bd)) \ar[d]^{H(\rho^x\ubul_Ap\ubul_Bd)} \ar@<.5ex>[lll]^-{\kappa}\\
H(x\bul_A(p\bul_Ad)) \ar@<.5ex>[d]^{H(\chi_{x\bul_Ap\bul_Bd})} \ar@<-.5ex>[d]_{H(\rho^x\ubul_Ap\ubul_Bd)} \ar@<.5ex>[rrr]^{H(\chi_{x\bul_Ap\bul_Bd})}
&&& H(x\bul_Ac\bul_A(p\bul_Bd)) \ar@<.5ex>[d]^{H(\chi_{x\bul_Ap\bul_Bd}\ubul_Bq\ubul_Ap\ubul_Bd)}  \ar@<-.5ex>[d]_{H(\rho^x\ubul_Ac\ubul_Ap\ubul_Bd)} \ar@<.5ex>[lll]^-{\xi_{x\bul_Ap\bul_Bd}}\\
H(x\bul_Ac\bul_A(p\bul_Bd)) \ar@<.5ex>[rrr]^-{H(\chi_{x\bul_Ac\bul_Ap\bul_Bd})}
&&& H(x\bul_Ac\bul_Ac\bul_A(p\bul_Bd)) \ar@<.5ex>[lll]^-{\xi_{x\bul_Ac\bul_Ap\bul_Bd}}
}
\end{equation}
Using the naturality of $\xi$, we find that the pair $(H(\chi_{m\bul_Ap\bul_Bd}), H(\rho^m\bul_Ap\bul_Bd))$ is equalized by $\xi_{m\bul_Ap\bul_Bd}\circ H(\rho^x\ubul_Ap\ubul_Bd)$. Hence, from the universal property of the equalizer, we find a morphism $\kappa:H(x\bul_A(p\bul_Bd)))\to H(x\bul^c(p\bul_Bd))$ such that $H({\rm eq}_x)\circ \kappa=\xi_{m\bul_Ap\bul_Bd}\circ H(\rho^x\ubul_Ap\ubul_Bd)$. Moreover,
we can compute
\begin{eqnarray*}
H({\rm eq}_x)\circ \kappa\circ H({\rm eq}_x) &=&  \xi_{m\bul_Ap\bul_Bd}\circ H(\rho^x\ubul_Ap\ubul_Bd)\circ H({\rm eq}_x)\\
&=& \xi_{m\bul_Ap\bul_Bd}\circ H(\chi_{m\bul_Ap\bul_Bd})\circ H({\rm eq}_x)= H({\rm eq}_x).
\end{eqnarray*}
Since the left column of \equref{equalizerdiagram} is an
equalizer, $H({\rm eq}_x)$ is a monomorphism, therefore $H({\rm
eq}_x)\circ \kappa$ is the identity. Hence, we find that the
equalizer in the left column is contractible by the maps
\[
\xymatrix{
H(x\bul_Ac\bul_Ap\bul_Bd) \ar[rr]^-{\xi_{m\bul_Ap\bul_Bd}} &&
H(x\bul_Ap\bul_Bd) \ar[rr]^-{\kappa} && H(x\bul^c(p\bul_Bd)).
}
\]
Since a contractible equalizer is preserved by any functor, this
equalizer is in particular preserved by $\qQ'$. Therefore, we find
that $\qQ' H= H' U_{\Omega,\cc} \qQ$ applied to the an equalizer
of the form \equref{cotensor} is an equalizer, that is
$H'((x\bul^c(p\bul_Bd))\bul_Bq)$ is the equalizer of the pair
$(H'(\rho^x\ubul_Ap\ubul_Bd\ubul_Bq),H'(\chi_{m\bul_Ap\bul_Bd}\bul_Bq))$.
Now consider the diagram
\begin{equation}\eqlabel{trivialeq}
\xymatrix{
x \ar[rr]^-{\rho^x} && x\bul_Ac \ar[d]^\cong \ar@<.5ex>[rr]^-{\rho^x\ubul_Ac} \ar@<-.5ex>[rr]_{x\ubul_A\Delta_c} && x\bul_Ac\bul_Ac \ar[d]^\cong\\
&& x\bul_Ap\bul_Bd\bul_Bq \ar@<.5ex>[rr]^-{\rho^x\ubul_Ap\ubul_Bd\ubul_Bq} \ar@<-.5ex>[rr]_-{\chi_{m\bul_Ap\bul_Bd}\bul_Bq} && x\bul_Ac\bul_Ap\bul_Bd\bul_Bq.
}
\end{equation}
This diagram defines an equalizer in $\Rcom(\Omega,\cc)$, which is
preserved by the forgetful functor $U_{\Omega,\cc}$. Moreover,
this equalizer is even split in $\Hom(\Omega,A)$, by the maps
\[
\xymatrix{
x && x\bul_Ac \ar[ll]_-{x\bul_A\varepsilon_c} && x\bul_Ac\bul_Ac. \ar[ll]_-{x\bul_Ac\bul_A\varepsilon_c}
}
\]
Therefore, this equalizer is preserved by the functor $H'$, and we
find in $\Yy$ that also $H'(x)$ is the equalizer of the pair
$(H'(\rho^x\ubul_Ap\ubul_Bd\ubul_Bq),H'(\chi_{m\bul_Ap\bul_Bd}\bul_Bq))$.
Therefore $H'(x)\cong H'((x\bul^c(p\bul_Bd))\bul_Bq)$. Since $H'$
reflects isomorphisms, we obtain that
$(x\bul^c(p\bul_Bd))\bul_Bq\cong x\cong x\bul^c(p\bul_Bd\bul_Bq)$
in $\Hom(\Omega,A)$, and therefore as well in $\Rcom(\Omega,\cc)$,
since the equalizer \equref{trivialeq} was preserved by
$U_{\Omega,\cc}$. Hence $\qQ$ preserves the equalizers of the form
\equref{cotensor}.

\ul{(2)}. By \thref{Qfullyfaithful} we only have to prove that the
equalizer $(d\bul_Bq)\bul^c(p\bul_Bd)$ is left $y$--pure for all
$y\in\Hom_1(\Omega,B)$. Since $\ul{\can}$ is an isomorphism, we can
identify $\cc$ with the comonad $(p\bul_Bd\bul_Bq)$, the same
apply for the comodules. Therefore, the equalizer
$(d\bul_Bq)\bul^c(p\bul_Bd)$ can be understood as the equalizer of
the pair
$(d\bul_Bq\bul_Ap\bul_B\gamma,\gamma\bul_Bq\bul_Ap\bul_Bd)$, where
$\gamma=(d\bul_B\eta\bul_Bd)\circ\Delta$. We have to prove that
the following diagram, where all subscripts have been removed by
typography's needs, is an equalizer in $\Rcom(\Omega,\dd)$.
\begin{equation}\eqlabel{equalizerydq}
\hspace{-1cm}\xymatrix{ \stackrel{\displaystyle y\bul((d\bul
q)\bul^c(p\bul d))} {} \ar[r]^-{y\bul\gamma} & y\bul(d\bul q)\bul
(p\bul d) \ar@<.5ex>[rr]^-{y\bul d\bul q\bul p\bul \gamma}
\ar@<-.5ex>[rr]_-{y\bul \gamma\bul q\bul p\bul d} && y\bul d\bul
q\bul p\bul d\bul q\bul p\bul d }
\end{equation}
Here we used that $\ol{\can}$ is an isomorphism. If we apply $H$
to \equref{equalizerydq}, then, taking into account that
$y\bul((d\bul q)\bul^c(p\bul d)) \cong y\bul d $, we obtain the
following diagram in $\Aa$.
\[
\xymatrix{ H(y\bul d)\ar[rr]^-{H(y\bul\gamma)} && H(y\bul d\bul
q\bul p\bul d) \ar@(ur,ul)[ll]_{\xi_{y\bul d}}
\ar@<.5ex>[rr]^-{H(y\bul d\bul q\bul p\bul\gamma)}
\ar@<-.5ex>[rr]_-{H(y\bul\gamma\bul q\bul p\bul d)}  && H(y\bul
d\bul q\bul p\bul d\bul q\bul p\bul d) \ar@(ur,ul)[ll]_{\xi_{y\bul
d\bul q\bul p\bul d}} }
\]
Let us show that this is a contractible equalizer, so in
particular an equalizer. The identities $\xi_{y\bul_Bd}\circ
H(y\bul_B\gamma)=H(y\bul_B\gamma)$ and
$\xi_{y\bul_Bd\bul_Bq\bul_Ap\bul_Bd}\circ
H(y\bul_Bd\bul_Bq\bul_Ap\bul_B\gamma)=H(y\bul_Bd\bul_Bq\bul_Ap\bul_Bd)$
follow directly from the $H$--separability of $(- \bul_B q)
U_{\Omega,\dd}$. Furthermore the naturality of $\xi$ implies that
$H(y\bul_B\gamma)\circ\xi_{y\bul_Bd}=\xi_{y\bul_Bd\bul_Bq\bul_Ap\bul_Bd}\circ
H(y\bul_B\gamma\bul_Bq\bul_Ap\bul_Bd)$.

Now compute the equalizer $(y\bul_Bd\bul_Bq)\bul^c(p\bul_Bd)$ of the pair
$(y\bul_Bd\bul_Bq\bul_Ap\bul_B\gamma,y\bul_B\gamma\bul_Bq\bul_Ap\bul_Bd)$. Since $H$ preserves this equalizer, we obtain that $H(y\bul_Bd)\cong H((y\bul_Bd\bul_Bq)\bul^c(p\bul_Bd))$ in $\Aa$. As $H$ reflects isomorphisms, we obtain $y\bul_Bd\cong (y\bul_Bd\bul_Bq)\bul^c(p\bul_Bd)$ in $\Rcom(\Omega,\dd)$, i.e. $y\bul_B-$ preserves the equalizer \equref{equalizerydq}, or $(d\bul_Bq)\bul^c(p\bul_Bd)$ is left $y$--pure.
\end{proof}

\begin{corollary}\colabel{strHsep}
Let $\cc=(A,c,\Delta_c,\varepsilon_c)$ and $\dd=(B,d,\Delta_d,\varepsilon_d)$ be comonads in $\Bb$ and $(p,q):\dd\to\cc$ be a comonad-morphism with adjunction such that $\Rcom(\Omega,\dd)$ satisfies the equalizer condition for $p$.
Suppose there exist functors $H''$, $H'$ and $\qQ'$ that make the following diagram commutative.
\[
\xymatrix{
\Rcom(\Omega,\dd) \ar[rr]^-\qQ
\ar[d]_{U_{\Omega,\dd}} && \Rcom(\Omega,\cc) \ar[d]^{U_{\Omega,\cc}} \\
\Hom(\Omega,B) \ar[rr]^-{- \bul_B q} \ar[d]_{H''} && \Hom(\Omega,A) \ar[d]^{H'}\\
\Bb \ar[rr]_{\qQ'} && \Yy
}
\]
\begin{enumerate}
\item The functor $\pP$ is fully faithful if the following conditions hold.
\begin{enumerate}[(i)]
\item $\ul{\can}$ is an isomorphism;
\item the functor $- \bul_B q$ is $H''$--separable;
\item $H''$ preserves equalizers of the form \equref{cotensor};
\item $H'$ reflects isomorphisms.
\end{enumerate}
\item $\qQ$ is fully faithful and therefore $(\qQ,\pP)$ is an equivalence of categories if in addition to (1) the following conditions hold,
\begin{enumerate}[(i)]
\item $\ol{\can}$ is an isomorphism (hence $(p,q)$ is a $\dd$-$\cc$ Galois comonad-morphism);
\item $H$ reflects isomorphisms.
\end{enumerate}
\end{enumerate}
\end{corollary}

\begin{proof}
\ul{(1)}. Since $\Rcom(\Omega,\dd)$ satisfies the equalizer
condition for $p$, $U_{\Omega,\dd}$ preserves preserves equalizers
of the form \equref{cotensor}. Put $H=H''\circ U_{\Omega,\dd}$,
then $H$ preserves equalizers of the form \equref{cotensor} as
well. By \leref{separable} (ii), we find that $(- \bul_B q)
U_{\Omega,\dd}$ is $H$--separable. The statement follows now from
\thref{Joyal} (1).\\
\ul{(2)} This follows in the same way from \thref{Joyal}(2).
\end{proof}

\section{Examples}\selabel{applications}

In this section we will briefly describe some situations of current
interest where our results apply.

\subsection{The bicategories of corings}\selabel{bicatcoring}

Consider $\Bb=\Bim(k)$ the bicategory of (unital) algebras over a
commutative ring, bimodules and homomorphisms of bimodules.
Comonads in $\Bim(k)$ are corings, and were studied, from the
point of view of bicategories, in \cite{BrzElkGT}. We will denote
the category $\Rcom(k,\cc)$ of right comodules over an $A$--coring
$\cc$ by $\Mm^\cc$, in analogy to the usual notation
$\Mm_A$ for the category of all right $A$--modules over a ring $A$.
An adjoint pair $(A,B,\Sigma,\Sigma^\dagger,\varepsilon,\eta)$ in
$\Bb$ was termed a comatrix coring context in \cite{BGT}. Since
$B$ is a ring with unit, we obtain that $\Sigma$ is finitely
generated and projective as a right $A$--module and $\Sigma^\dagger
\cong \Sigma^*$. By \equref{adjointmonad}, we have the  $A$--coring
$\Sigma^*\ot_B\dd\ot_B\Sigma$. This construction was considered in
\cite[Theorem 3.1]{BrzElkGT} and it generalizes that of finite
comatrix corings \cite{EGT:comatrix}.

If we apply the techniques developed in the previous sections to
the present situation, then we recover the situation studied in
\cite{BrzElkGT}, and, in particular, those in \cite[Section
5]{BrzElkGT} concerning the pull-back and push-out functors. From
\prref{adjuntos},
 \thref{weakstructure}, \thref{Qfullyfaithful} and
\thref{strongstructure} we deduce:

\begin{theorem}
Let $A$ and $B$ be unital rings, $\dd$ a $B$--coring that is flat
as a left $B$--module and $\cc$ an $A$--coring. Consider
$\Sigma\in{_B\Mm^\cc}$ such that $\Sigma$ is finitely generated
and projective as a right $A$--module, and let $\{ (e_i, f_i) \}
\subset \Sigma \times \Sigma^*$ be a finite dual basis.
\begin{enumerate}[(1)]
\item We have a pair of adjoint functors $(F,G)$
\begin{eqnarray*}
F: &\Mm^\dd \to \Mm^\cc,& F(M)=M\ot_B\Sigma; \\
G: &\Mm^\cc \to \Mm^\dd,& G(N)=N\ot^\cc(\Sigma^*\ot_B\dd).
\end{eqnarray*}
\item The functor $G$ is fully faithful if and only if the
canonical map
$$\can : \Sigma^*\ot_B\dd\ot_B\Sigma\to \cc, \qquad
\can(\varphi\ot_Bd\ot_Bu)=\varphi(\varepsilon_\dd(d)u_{[0]})u_{[1]}$$
is an isomorphism and
$$(N\ot^\cc(\Sigma^*\ot_B\dd))\ot_B\Sigma\cong
N\ot^\cc(\Sigma^*\ot_B\dd\ot_B\Sigma)$$ for all $N\in\Mm^\cc$.
\item The functor $F$ is fully faithful if and only if the map
$$\ol{\can}:\dd\to (\dd\ot_B\Sigma)\ot^\cc(\Sigma^*\ot_B\dd),
\qquad \ol{\can}(d)=d_{(1)}\ot_Be_i\ot_Af_i\ot_Bd_{(2)}$$ is an
isomorphism and the map
$\dd\to(\dd\ot_B\Sigma)\ot_A(\Sigma^*\ot_B\dd)$ is a pure morphism
in ${_B\Mm}$. \item $(F,G)$ is an equivalence of categories if and
only if $\can$ is an isomorphism, $-\otimes_B\Sigma:\Mm^\dd\to
\Mm^\cc$ reflects isomorphisms and
$$(N\ot^\cc(\Sigma^*\ot_B\dd))\ot_B\Sigma\cong
N\ot^\cc(\Sigma^*\ot_B\dd\ot_B\Sigma)$$ for all $N\in\Mm^\cc$.
\end{enumerate}
\end{theorem}

\begin{remark}
Item (4) in the previous theorem has the following alternative
formulation \cite[Theorem 5.2]{BrzElkGT}: $(F,G)$ is an
equivalence of categories and $\cc$ is flat as a left $A$--module
if and only if $\can$ is an isomorphism and $\dd\ot_B\Sigma$ is
faithfuly coflat as a left $\dd$--comodule.~\qed
\end{remark}

\subsection{Coendomorphism corings and Morita-Takeuchi Theory}\selabel{coend}

Let us briefly argue how Morita-Takeuchi Theory on equivalences of
categories of comodules over corings (see \cite{Zarouali:phd} and
\cite[Section 23]{BrzWis:book}) might be derived from our set up.
Let $\cc$ (resp.\ $\dd$) be a coring over a $k$--algebra $A$ (resp.\
$B$). It follows from \cite[Proposition 3.4]{Gomez:2002} (see e.g.\
\cite[Theorem 2.1.13]{Zarouali:phd}) that for corings flat as left
modules over their ground rings, any equivalence between their
categories of right comodules is given by a cotensor product
functor. Thus, let us consider bicomodules ${}_\cc P_\dd$ and
${}_\dd Q_\cc$, and assume that $ - \otimes^\dd Q : \Mm^\dd
\rightarrow \Mm^\cc$ is left adjoint to $- \otimes^\cc P : \Mm^\cc
\rightarrow \Mm^\dd$. If ${}_B\dd$ is flat, then, by
\cite[Proposition 4.2]{Gomez:2002}, this is equivalent to assume
that $P$ is $(A, \dd)$--quasi-finite, that is, $- \otimes^\dd Q:
\Mm^\dd \rightarrow \Mm_A$ is left adjoint to the functor $-
\otimes_A P: \Mm_A \rightarrow \Mm^\dd$. This situation is
modelled by our theory in the framework of the $2$--category
$\CAT$ of categories, functors and natural transformations as
follows. Let $*$ be the discrete one-object category. Consider the
category $\Mm^\dd$ as a $0$--cell, and the trivial comonad on
$\Mm^\dd$ (built on the identity functor $id_{\Mm^\dd}: \Mm^\dd
\rightarrow \Mm^\dd$). Then the category $\Rcom(*,id_{\Mm^{\dd}})$
is
 isomorphic to $\Mm^{\dd}$. By considering the comonad $-\otimes_A
 \cc$ on the category $\Mm_A$, we get an isomorphism of categories
 $\Rcom(*,-\otimes_A \cc) \simeq \Mm^\cc$. In this way, we have
 an adjoint pair $(p,q)$ in $\CAT$, where $p = - \otimes_A P$, $q
 = - \otimes^{\dd} Q$ (the horizontal composition in $\CAT$ is the
 opposite of the composition of functors, thinking that they act on the left on objects and
 morphisms). Moreover, since $P$ is a $\cc$-$\dd$--bicomodule, we
 obtain that $(p,q)$ is a comonad-morphism from $id_{\Mm^\dd}$ to
 $- \otimes_A \cc$, which gives, by \prref{locdual} (see also \cite[Proposition 2.3]{GT:comonad}), a
 homomorphism of comonads $(-\otimes_A P) \otimes^\dd Q
 \rightarrow - \otimes_A \cc$. Since $- \otimes^\dd Q$ is a left
 adjoint, it preserves equalizers, whence, by \cite[Lemma
 2.2]{Gomez:2002}, there is a natural isomorphism $(-\otimes_A P) \otimes^\dd
 Q \simeq - \otimes_A ( P \otimes^\dd Q)$. In this way, the
 foregoing homomorphism of comonads is determined by a
 homomorphism of $A$--corings $P \otimes^\dd Q \rightarrow \cc$.
 Now, we can easily deduce from \thref{strongstructure} or \cite[Theorem 2.7]{GT:comonad}
 that, if ${}_B\dd$ is flat, then the  functors $- \otimes^\cc P :
 \Mm^\cc \rightarrow \Mm^\dd$ and $- \otimes^{\dd} Q : \Mm^\dd
 \rightarrow \Mm^\cc$ give an equivalence of categories if and only
 if $- \otimes^\dd Q$ is exact and faithful. From this, it is easy
 to deduce \cite[23.12]{BrzWis:book}.

\begin{remark}
Both the examples given in \seref{bicatcoring} and \seref{coend}
describe equivalences between two categories of comodules. The
relation between the two theories is clarified as follows.
Starting with an $A$--coring $\cc$ and a $B$--coring $\dd$, and a
comatrix coring context
$(A,B,\Sigma,\Sigma^\dagger,\varepsilon,\eta)$, the comodule
$Q=\Sigma^\dagger\ot_B\dd$ is $(A,\dd)$--quasi-finite and
$-\ot_AQ:\Mm_A\to\Mm^\dd$ has a left adjoint, represented by the
$\dd$-$A$--bicomodule $P=\dd\ot_B\Sigma$. In particular, we find
the following isomorphism between the two associated comatrix
corings occuring in both theories $Q\ot^\dd
P=(\Sigma^\dagger\ot_B\dd)\ot^\dd(\dd\ot_B\Sigma)\cong
\Sigma^\dagger\ot_B\dd\ot_B\Sigma$.~\qed
\end{remark}

\subsection{Comatrix Corings over firm rings}\selabel{firm}

A (non-unital) associative ring $R$ is called firm if the
multiplication on $R$ induces an isomorphism $R\ot_RR\cong R$. A
right $R$--module over a firm ring $R$ is called firm if and only
if $M\ot_RR\cong R$. One can easily construct a bicategory
$\Frm(k)$, whose $0$--cells are firm rings, $1$--cells are firm
bimodules and $2$--cells are homomorphisms of bimodules. It was
proven in \cite{Ver:equiv} that firm rings can be characterized as
corings over their Dorroh-extension. The Galois theory for corings
over firm rings has been initiated in \cite{GTV},
\cite{GT:comonad} and, in a more profound treatment, in
\cite{V:PhD}. The situation studied in \cite{GTV} is subsumed by
the theory developed in the present paper by taking $\Bb=\Frm(k)$,
and the main results of \cite{GTV} are obtained as consequences.
The version of \cite{BrzElkGT} for corings over firm rings is also
recovered in this way. We leave the details of these constructions
to the reader.

\subsection{Descent theory}\selabel{descent}
It was pointed out in \cite[Sect. 4.8]{CMZ} and \cite[Example
2.1]{Brz:structure} that the category of descent data associated to
a ring extension $B \rightarrow A$ is isomorphic to the category of
right comodules over the $A$--coring $A \otimes_B A$. As mentioned
in the introduction, some authors would prefer to name all the
theory we have developed descent theory, however we reserve this
name for the following special situation. Let $A$ and $B$ be
$k$--algebras (with unit) and $\Sigma$ a $B$-$A$--bimodule (in the
case of the ring extension $B \rightarrow A$, take $\Sigma =
{}_BA_A$). Then we can consider the functor
$$-\ot_B\Sigma:\Mm_B\to\Mm_A.$$
The descent problem in this setting is described as the following
question
\begin{quote}
Which right $A$--modules are of the form $N\ot_B\Sigma$ for some
$N\in\Mm_B$, i.e. for which $M\in\Mm_A$, can we find an $N\in\Mm_B$
such that $M\cong N\ot_B\Sigma$ ?
\end{quote}
A solution to the problem can be formulated if $\Sigma$ is finitely
generated and projective as a right $A$--module; then we can
construct the comatrix $A$--coring $\cc=\Sigma^*\ot_B\Sigma$ and
consider its category of comodules $\Mm^\cc$ as the category of
descent data (see \cite[\S 5.2]{BrzElkGT} for the notion of a
generalized descent datum). The functor $F$ factorizes in the
following way
\[
\xymatrix{
F:\Mm_B\ar[r]^-{-\ot_B\Sigma} & \Mm^\cc\ar[r]^-{U_\cc} & \Mm_A,
}
\]
where $U_\cc$ denotes as usual the forgetful functor.  Moreover,
$-\ot_B\Sigma$ establishes an equivalence between $\Mm_B$ and
$\Mm^\cc$ if $\Sigma$ is faithfully flat as a left $B$--module
\cite[Theorem 3.10]{EGT:comatrix}. Conversely, if $\Sigma^*
\otimes_B \Sigma$ is flat as a left $A$--module (e.g., if
${}_B\Sigma$ is flat), and $-\ot_B\Sigma$ induces an equivalence
between $\Mm_B$ and $\Mm^\cc$, then $\Sigma$ is faithfully flat as a
left $B$--module \cite[Theorem3.10]{EGT:comatrix}. The ``classical''
descent theorem for a noncommutative ring extension was proved in
\cite[Teorema 8]{Cip} (see also \cite{Nuss} and \cite[Theorem
3.8]{Brz:structure}).

\thref{Joyal} can be applied in this situation. A remarkable fact is that in under certain conditions, part (2) of \thref{Joyal} has a converse.
Let us first prove the following lemma.

\begin{lemma}\lelabel{Hseparabledescent}
Let $A$ and $B$ be rings and $\Sigma\in{_B\Mm_A}$ finitely generated
and projective as right $A$--module with finite dual basis
$e=e_i\ot_Af_i\in \Sigma\ot_A\Sigma^*$. Consider
$H=\Hom_\ZZ(-,\QQ/\ZZ):\Mm_B\to {_B\Mm}^{op}$. The following
statements are equivalent:
\begin{enumerate}[(i)]
\item The functor $-\ot_B\Sigma:\Mm_B\to\Mm_A$ is $H$--separable;
\item the map $H(\eta):H(\Sigma\ot_A\Sigma^*)\to H(B)$ is split epi in ${_B\Mm_B}$, where $\eta(b)=eb=be$.
\end{enumerate}
\end{lemma}

\begin{proof}
We only have to prove that (ii) implies (i). We have to check that
for all $M\in\Mm_B$ the following morphism has a right inverse
$$H(\eta_M):H(M\ot_B\Sigma\ot_A\Sigma^*)\to H(M),$$
where $\eta_M(m)=m\ot_Be$. Using natural isomorphisms (in the
vertical arrows), we see in the following diagram that the right
inverse $\xi$ of $\eta$ induces a splitting map $\xi_M$ for
$\eta_M$.
\[
\xymatrix{
H(M\ot_B\Sigma\ot_A\Sigma^*)\ar@<.5ex>[rr]^-{H(\eta_M)}
 \ar[d]_-{\cong} && H(M) \ar[d]^-{\cong}\ar@<.5ex>[ll]^-{\xi_{M}} \\
\Hom_B(M,H(\Sigma\ot_A\Sigma^*)) \ar@<.5ex>[rr]^{H(\eta)\circ -}
 && \Hom_B(M,H(B)) \ar@<.5ex>[ll]^{\xi\circ-}
}\vspace{-32pt}{~}
\]
\end{proof}
\vspace{7pt}

\begin{theorem}[{\cite[Theorem 2.7]{CDV}, \cite[Theorem 4.1]{Mes:extension}}]\thlabel{CDV}
Let $A$ and $B$ be rings, take $\Sigma\in{_B\Mm_A}$ finitely generated and projective as right $A$--module.
Denote by $e\in\Sigma\ot_A\Sigma^*$ the finite dual basis for $\Sigma$ and put $\cc=\Sigma^*\ot_B\Sigma$. Then $(i)$ implies $(ii)$ implies $(iii)\Leftrightarrow(iv)$.
\begin{enumerate}[(i)]
\item $H(\eta):H(\Sigma\ot_A\Sigma^*)\to H(B)$ is split epi in ${_B\Mm_B}$;
\item $-\ot_B\Sigma:\Mm_B\to\Mm^\cc$ is an equivalence of categories;
\item $-\ot_B\Sigma:\Mm_B\to\Mm^\cc$ is fully faithful;
\item $\eta:B\to\Sigma\ot_A\Sigma^*$, $\eta(b)=be=eb$ is pure in ${_B\Mm}$.
\end{enumerate}
If $B$ is commutative and every $\varphi\in{\End_A}(\Sigma)$ is left
$B$--linear, or $B$ is a separable $k$--algebra, then $(iv)$ implies
$(i)$ and all statements become equivalent.
\end{theorem}
\begin{proof}
$\ul{(i)\Rightarrow(ii)}$.
From \leref{Hseparabledescent}, we know that $-\ot_B\Sigma$ is $H$--separable. Moreover, it is known that the functor $H$ preserves all equalizers
and reflects isomorphisms, as $\QQ/\ZZ$ is an injective cogenerator in the category $\Ab$ of abelian groups. Hence, the implication follows by \thref{Joyal}, taking $\dd=B$, the trivial comonad, $q=\Sigma$,  $p=\Sigma^*$, $\cc=\Sigma\ot_A\Sigma^*$, the comatrix coring associated to $\Sigma$,  $\Aa=\Bb={_B\Mm}^{op}$ and $\Yy=\Ab^{op}$, $H=\Hom_\ZZ(-,\QQ/\ZZ):\Mm_B\to {_B\Mm}^{op}$, $H'=\Hom_\ZZ(-,\QQ/\ZZ):\Mm_A\to \Ab^{op}$ and $\qQ'=\Hom_B(\Sigma,-): {_B\Mm}^{op}\to\Ab^{op}$.\\
$\ul{(ii)\Rightarrow(iii)}$. Trivial.\\
$\ul{(iii)\Leftrightarrow(iv)}$. By \thref{weakstructure}.\\
$\ul{(iv)\Rightarrow(i)}$ This last implication is proven in \cite[Proposition 2.6]{CDV} in the case $B$ is commutative, and in \cite[Theorem 3.4 $(i)\Rightarrow(iii)$]{Mes:extension} in the more general case $B$ is a separable $k$--algebra.
\end{proof}

If $B\to A$ is an extension of commutative rings and we take $\Sigma=A$, then \thref{CDV} reduces to the Joyal-Tierney theorem.

\subsection{The Hopf-Galois theories}

The Galois theory for corings, originated, with many intermediate
steps, from Hopf-Galois theory, which itself initiated as a
generalization of classical Galois theory. The connection between
corings and entwining structures was made explicit in
\cite{Brz:structure}, where it was shown that the category of
(right) entwined modules for an algebra $A$ entwined with a
coalgebra $C$ is isomorphic to the category of right comodules
over a suitable $A$--coring built on $A \otimes C$. This in
particular covers the case of a comodule algebra $A$ over the Hopf
algebra $H$, which is an $H$--Galois extension of $A^{{\sf
co}H}=\{a\in A~|~\rho(a)=a\ot 1_H\}$ if and only if the canonical
map
$$\can:A\otimes_{A^{{\sf co}H}}A\to A\otimes_k H,\quad\can(a\otimes a')=aa'_{[0]}\otimes a'_{[1]},$$
is an isomorphism, where $\rho^A(a)=a_{[0]}\ot a_{[1]}$ denotes the
$H$--coaction on $A$. This canonical map is nothing but the canonical
map corresponding to the grouplike element $1 \otimes 1$ of the
$A$--coring $A \otimes H$. Therefore, Hopf-Galois theory can be
generalized in the framework of corings. This theory can be extended
in many ways. One can replace $H$ by a weak Hopf algebra
\cite{CDG:galweak}, a Hopf algebroid \cite{Bohm:gal} or a coalgebra
\cite{BrzezinskiH}. We refer to e.g. \cite{Cae:descent} and
\cite{Wis:galfields} for more detailed overviews.

\subsection{Comatrix Corings over Quasi-Algebras}

In this section we provide a new type of examples where our general
theory applies, associated to (dual) quasi-bialgebras. The
philosophy is similar to the construction of comatrix corings over
firm rings \cite{GTV}, where we constructed comatrix corings,
replacing one of the unital rings by a firm ring. Firmness means
exactly that the category of bimodules over the firm ring is a
monoidal category with the base ring as monoidal unit. In the same
way, we can work with a ring (with unit), which is not associative,
but inducing a canonical isomorphism
$M\ot_R(N\ot_RP)\cong(M\ot_RN)\ot_RP$ for all $R$--bimodules $M$, $N$
and $P$. Such a framework can be obtained in the following setting.

Let $k$ be a commutative ring, unadorned tensor products in this
section are tensor products over $k$. A dual quasi-bialgebra is a
sextet $(H,\Delta,\epsilon,\mu,\eta,\phi)$, where $\Delta:H\to H\ot
H$ is a coassociative coproduct and $\varepsilon:H\to k$ is a counit
for $\Delta$. Furthermore $\phi:H\ot H\ot H\to k$ is a unital
$3$--cocycle that is convolution invertible, this means that the
following identities hold for all $a,b,c,d\in H$
\begin{eqnarray}
&&\hspace{-4cm}\phi(b_{(1)},c_{(1)},d_{(1)})\phi(a_{(1)},b_{(2)}\cdot c_{(2)},d_{(2)})\phi(a_{(2)},b_{(3)},c_{(3)})\eqlabel{cocycle1}\\
&=& \phi(a_{(1)},b_{(1)},c_{(1)}\cdot d_{(1)})\phi(a_{(2)}\cdot b_{(2)},c_{(2)},d_{(2)});\nonumber\\
\phi(a,1,b)&=&\varepsilon(a)\varepsilon(b).\eqlabel{cocycle2}
\end{eqnarray}
There exists a map $\phi^{-1}:H\ot H\ot H\to k$ such that
\begin{eqnarray}
\phi(a_{(1)},b_{(1)},c_{(1)})\phi^{-1}(a_{(2)},b_{(2)},c_{(2)})&=&\varepsilon(a)\varepsilon(b)\varepsilon(c)\eqlabel{convinv}\\
&=&\phi^{-1}(a_{(1)},b_{(1)},c_{(1)})\phi(a_{(2)},b_{(2)},c_{(2)}).\nonumber
\end{eqnarray}
Furhtermore, the product $\mu:H\ot H\to H,\mu(a\ot b)=a\cdot b$ is
associative up to conjugation with by $\phi$, i.e.
$$a_{(1)}\cdot(b_{(1)}\cdot c_{(2)})\phi(a_{(2)},b_{(2)},c_{(3)})
=\phi(a_{(1)},b_{(1)},c_{(2)}) (a_{(2)}\cdot b_{(2)})\cdot c_{(3)}$$
for all $a,b,c\in H$.

A right comodule over a dual quasi-bialgebra, is a right comodule
over $H$, i.e. a $k$--module $M$ together with a map $\rho^M:M\to
M\ot H$ satisfying $(M\ot \Delta)\circ\rho^M=(\rho^M\ot
H)\circ\rho^M$ and $(M\ot\varepsilon)\circ\rho^M$.

Recall \cite{Bohm:int} that, in analogy to $\Bim(k)$, we can
construct out of any monoidal category with coequalizers $\Mm$, a
bicategory $\Bim(\Mm)$, whose $0$--cells are algebras in $\Mm$,
$1$--cells are bimodules over these algebras and $2$--cells are
bilinear maps. It is known that the category $\Mm^H$ of right
comodules over a dual quasi bialgebra form a monoidal category. The
tensor product of $\Mm^H$ is $\ot=\ot_k$, the monoidal unit is $k$,
whose coaction is given by $\eta:k\to H$ and the associativity
constraint in $\Mm^H$ is given by the following natural isomorphism,
$$\Phi_{M,N,P}:(M\ot N)\ot P\to M\ot(N\ot P), \Phi((m\ot n)\ot p)=m_{[0]}\ot(n_{[0]}\ot p_{[0]})\phi(m_{[1]},n_{[1]},p_{[1]}).$$
Since $\Mm^H$ is a comodule category, it has coequalizers, therefore, we can construct the bicategory $\Bim(\Mm^H)$. Moreover, we can apply Galois theory as developed in earlier sections to this bicategory, which will be discussed in a forthcoming paper.
Let us just state the following remarkable result, that allows the explicit construction of a new type of comatrix corings.

Consider an associative $k$--algebra $A$, an $H$--comodule
quasi-algebra $B$ (i.e. $B$ is an algebra in the monoidal category
$\Mm^H$) and let $(A,B,\Sigma,\Sigma',\mu,\eta)$ be an adjoint
pair in $\Bim(\Mm^H)$. This means that $\Sigma\in{_B\Mm^H_A}$,
$\Sigma'\in {_A\Mm^H_B}$, $\eta:B\to \Sigma\ot_A\Sigma'$ is a
morphism in ${_B\Mm^H_B}$ and $\mu:\Sigma'\ot\Sigma\to A$ is a
morphism in ${_A\Mm^H_A}$. Since the map $\eta$ is completely
determined by its image on $1_B$, let us denote
$\eta(1_B)=e_i\ot_Af_i\in \Sigma\ot_A\Sigma'$. The conditions of
an adjoint pair then translate into the following identities
\begin{equation}\eqlabel{quasicc}
\begin{array}{rcl}
x&=&e_{i[0]}\mu(f_{i[0]}\ot_B x_{[0]})\phi(e_{i[1]},f_{i[1]},x_{[1]});\\
y&=&\mu(y_{[0]}\ot_Be_{i[0]})f_{i[0]}\phi^{-1}(e_{i[1]},f_{i[1]},x_{[1]}),
\end{array}
\end{equation}
for all $x\in\Sigma$ and $y\in\Sigma'$. By \equref{adjointmonad}
we immediately obtain the following result, which generalizes the
construction of (finite) comatrix corings given in
\cite[Proposition 2.1]{EGT:comatrix}.

\begin{theorem}
With notation introduced above,
$\Sigma'\ot_A\Sigma$ is an $A$--coring, with counit $\mu$ and comultiplication $\Delta$ defined for all $y\ot_Bx\in\Sigma'\ot_B\Sigma$ by
\begin{eqnarray*}
\Delta(y\ot_B x)&=&(y_{[0]}\ot_B e_{i[0]})\ot_A(f_{i[0]}\ot_B x_{[0]})\phi(y_{[1]}\cdot e_{i[1]},f_{i[1]},x_{[1]})\phi^{-1}(y_{[2]},e_{i[2]},f_{i[2]})\\
&=& (y_{[0]}\ot_B e_{i[0]})\ot_A(f_{i[0]}\ot_B x_{[0]})
\phi^{-1}(y_{[1]},e_{i[1]},f_{i[1]}\cdot x_{[1]})
\phi(e_{i[2]},f_{i[2]},x_{[2]}).
\end{eqnarray*}
\end{theorem}

\subsection{Galois theory of Group-corings}

In \cite{CJW} a Galois theory for group corings is being developed.
We will now show that this theory fits as well in our general
framework. To this end, we introduce a new bicategory, that unifies
the bicategory of bimodules with the notion of a Turaev category,
introduced in \cite{CDL}.

\begin{Definition}
Let $k$ be a commutative ring, then we define the bicategory $\Tur(k)$ as follows. The $0$--cells are $k$--algebras, a $1$--cell $\ul{M}:A\to B$ is a couple $(X,(M_x)_{x\in X})$, where $X$ is a set and all $M_x$ are $A$-$B$--bimodules. A $2$--cell $\ul{\varphi}=(f,\varphi):\ul{M}=(X,M_x)\to\ul{N}=(Y,N_y)$ consists of a map $f:Y\to X$ together with a collection of $A$-$B$--bilinear maps $\varphi_y:M_{f(y)}\to N_y$. The composition of $1$--cells is defined as follows. Take $\ul{M}=(X,M_x)\in\Hom_1(A,B)$ and $\ul{N}=(Y,N_y)\in\Hom_1(B,C)$ then we define $\ul{M}\bul_B\ul{N}=(X\times Y,M_x\ot_BN_y)$. In the same way, we define the horizontal composition $\ubul$ of $2$--cells. For $\ul{\varphi}=(f,\varphi):\ul{M}\to\ul{N}$ and $\ul{\psi}=(g,\psi):\ul{N}\to\ul{P}=(Z,P_z)$ we define the vertical composition $\ul{\psi}\circ\ul{\varphi}=(f\circ g, (\psi_z\circ\varphi_{g(z)})_{z\in Z})$.

If we consider $\Tur(k)$ with only one $0$--cell $k$, then we recover the notion of a Turaev category $\Tt_k$ of \cite{CDL}.
\end{Definition}

It is possible to construct another bicategory, with the same
$0$--cells and $1$--cells, but where $2$ cells are of the form
$\ul{\varphi}=(f,\varphi):\ul{M}=(X,M_x)\to\ul{N}=(Y,N_y)$
consisting of a map $f:X\to Y$ together with a collection of
$A$-$B$--bilinear maps $\varphi_x:M_{x}\to N_{f(y)}$. Considering such a
bicategory with only one $0$--cell $k$, we recover the notion of a
Zunino category $\Zz_k$ introduced in \cite{CDL}.

Both bicategories introduced above, admit a locally faithful
pseudo functor from $\Bim(k)$. More precisely,
$$F:\Bim(k)\to \Tur(k)$$
is defined by $F(A)=A$ on $0$--cells, $F(M)=(\{*\},M)$ on $1$--cells.

As it was shown in \cite{CDL}, the Turaev and Zunino categories
lead to a conceptual interpretation of group-algebras, -coalgebras
and -Hopf algebras, being algebras, coalgebras and Hopf algebras
in a (braided) monoidal category. In a similar way, we can use the
Turaev and Zunino bicategories to interpret other constructions of
`group'-type. For example, mixed distributive laws (entwining
structures) can be constructed in any bicategory. Their
interpretation in the Turaev bicategory lead to the notion of
`group entwining structures' as introduced in
\cite{Wang:groupentw}.

A comonad in $\Tur(k)$ coincides with the notion of a group-coring as defined in the recent paper \cite{CJW}.
If we apply the Galois theory of this paper to the Turaev bicategory, we obtain the Galois theory of \cite{CJW} as a special situation.

\subsection{Comonads over CAT}

Consider now $\Bb=\CAT$, a category whose $0$--cells are categories
which are small in some Grothendieck universe, the $1$--cells are
functors, and the $2$--cells are natural transformations between
them. Then comonads are cotriples, and the Galois theory that we
have developed is tightly linked to the theory of comonadicity (or
cotripleability) of functors. In particular, we recover the famous
theorem of Beck. Let $(F,G)$ be a pair of adjoint functors with
$F:\Bb\to \Aa$ and $G:\Aa\to \Bb$. Recall that we can construct a
cotriple (comonad) $C=FG$ on $\Aa$ which induces a pair of adjoint
functors $(\Ff^C,\Gg^C)$ with $\Ff^C:\Aa^C\to\Aa$ and
$\Gg^C:\Aa\to\Aa^C$. Moreover there exists a unique functor
$K:\Bb\to\Aa^C$ such that $F=\Ff^CK$ and $KG=\Gg^C$. From
\thref{weakstructure} and \thref{strongstructure}, we immediately
obtain the following.
\begin{theorem}[Beck]
If the category $\Bb$ has equalizers, then $K$ has a left adjoint $L:\Aa^C\to \Bb$, which is fully faithful if and only if $F$ preserves equalizers. If moreover $F$ reflects isomorphisms, then $K$ is an equivalence between the categories $\Aa^C$ and $\Bb$.
\end{theorem}

This part of the theory and its connection with the theory of Galois comodules over firm rings has been discussed in more detail in \cite{GT:comonad}.

Let us just remark that Galois theory in $\CAT$ applied to the
situation of corings and comodules, recovers the theory of Galois
comodules in the sense of Wisbauer \cite{Wis:galcom}, termed
comonadic-Galois comodules in \cite{Ver:equiv}. Consider any firm
ring $R$ and an $A$--coring $\cc$ and take $\Sigma\in{_R\Mm^\cc}$.
Then we say that $\Sigma$ is an $R$-$\cc$ comonadic-Galois comodule
if the following morphism is an isomorphism for all $M\in\Mm_A$
\begin{equation}\eqlabel{canWis}
\can_M:\Hom_A(\Sigma,M)\ot_R\Sigma\to M\ot_A\cc,\quad\can_M(\varphi\ot_Ru)=\varphi(u_{[0]})\ot_Au_{[1]}.
\end{equation}
For any $R$-$\cc$--bicomodule $\Sigma$, the functor
$-\ot_R\Sigma:\Mm_R\to\Mm_A$ has a right adjoint given by
$\Hom_A(\Sigma,-)\ot_RR:\Mm_A\to\Mm_R$. Hence we can construct the
associated comonad in $\CAT$ and compare it with $\cc$ by a
canonical comonad-morphism, which becomes exactly the canonical
cotriple morphism \equref{canWis}. To apply the Galois theory in
$\Frm(k)$ and obtain a firm Galois comodule, we need however a
comonad-morphism with adjunction in $\Frm(k)$, which means that
there must exist a $2$--cell in $\Frm(k)$ that represents the
functor $\Hom_A(\Sigma,-)\ot_RR$. By the Eilenberg-Watt's theorem
over firm rings (see \cite[Theorem 3.1]{Ver:equiv}), we know that
this means exactly that $\Sigma$ is $R$--firmly projective as right
$A$--module. This condition is satisfied once
$\Hom_A(\Sigma,-)\ot_RR:\Mm_A\to\Mm_R$ has a right adjoint, so in
particular if $\Hom^\cc(\Sigma,-)\ot_RR:\Mm^\cc\to\Mm_R$ has a
right adjoint or if $-\ot_R\Sigma:\Mm_R\to\Mm^\cc$ is an
equivalence of categories. In these situations, both theories
coincide. A consequence of this reasoning is that any equivalence
of categories between a category of comodules and a category of
modules over a unital ring reduces to `finite' Galois theory (i.e.
in the sense of \cite{EGT:comatrix}), this was proven in
\cite{Ver:equiv}.

\subsection{Galois theory for monads}
Since monads in a bicategory $\Bb$ are comonads in the bicategory $\Bb^{co}$, we obtain by direct dualization a Galois theory for monads. The explicit definition of the associated bicategories of Eilenberg-Moore objects have been given in \cite{LacStr} and were one of our major inspirations. However, the Galois theory was not developed there and can therefore be derived from our approach.

An interesting example of this dual theory can be the theory of
monadicity of functors, taking $\Bb=\CAT$. Among many, an
interesting paper is \cite{JanTho:descIII}, where a categorical
interpretation of the Joyal-Tierney theorem is proven. However,
comparing their result to \thref{Joyal}, they only prove a version
of part (1) of the theorem.

\subsection{Galois theory of Matrix \emph{C}--rings}
Let $k$ be a field and consider the category $\Bic(k)$. An adjoint
pair in this bicategory consists of two $k$--coalgebras $C$ and
$D$, a $C$-$D$--bicomodule $M$, a $D$-$C$--bicomodule $N$ and two
bicolinear maps $\sigma:C\to N\ot^D M$ and $\tau:M\ot^CN\to D$,
satisfying $\tau\ot^DN\cong N\ot^C\sigma$ and $M\ot^C\sigma\cong
\tau\ot^DM$. Applying \equref{adjointmonad}, we find that
$M\ot^DN$ is a monad in $\Bb$ and $N\ot^CM$ is a comonad in $\Bb$.
If we consider the comonad, then we can apply our general theory
(or in fact, the special case of the coendomorphism corings), if
we consider the monad, then we have to apply the dual version of
the theory. This dual Galois theory is recently developed in
\cite{BrzTur:Crings}, the monad $M\ot^DN$ is termed a \emph{matrix
$C$--ring}.

\end{document}